\documentclass{amsart}

\usepackage{amsbsy,bm,fixmath}
\usepackage{epsfig}
\usepackage{amsfonts}
\usepackage{amssymb}
\usepackage{amsmath}
\usepackage{stmaryrd}
\usepackage{float}
\usepackage{xcolor}
\usepackage{comment}
\usepackage{dsfont}
\usepackage{verbatim}
\usepackage{mathtools}

\newtheorem{theorem}{Theorem}[section]
\newtheorem{lemma}[theorem]{Lemma}
\theoremstyle{definition}
\newtheorem{definition}[theorem]{Definition}
\newtheorem{example}[theorem]{Example}

\newtheorem{corollary}[theorem]{Corollary}
\theoremstyle{remark}
\newtheorem{remark}[theorem]{Remark}

\def\EE{\mathbb E}
\def\cal{\mathcal}
\numberwithin{equation}{section}

\begin{document}

\title[Mean field games]{Mean Field Games: Large Sparse Network Limits and  Laplexion Dynamics}\thanks{This work was supported by  NSERC (Grant RGPIN 2019-06171, M.H.),  AFOSR (Grant FA9550-23-1-0015, P.E.C.), and NSERC (Grant RGPIN 2019-05336, P.E.C.).}


\author{Minyi Huang}
\address{School of
Mathematics and Statistics, Carleton University,
Ottawa, ON, Canada, and GERAD, Montreal.  }
\email{mhuang@math.carleton.ca}

\author{Peter E. Caines}
\address{Department
of Electrical and Computer Engineering, McGill University, Montreal, QC, Canada, and GERAD, Montreal.  }
\email{peter.caines@mcgill.ca}

%    General info
%\subjclass[2020]{Primary xx, xx; Secondary xx, xx}

\date{}


\begin{abstract}
Dynamic games are considered with large subpopulations distributed
over large sparse graphs, where each agent has mean field coupling with all agents  within the same cluster and meanwhile
receives impact from  neighboring clusters. Tractable
limit models are obtained via the network limits of graph Laplacian
operators modelling second order interactions. The
resulting mean field game (MFG) has what is called Laplexion dynamics. For ring and torus  topology cases, solutions for linear quadratic Laplexion MFG systems with infinite populations on infinite node networks are obtained via a second order linear partial differential equation system. This work focusses on the relationship between the finite population model and the sparse network limit model.  We  prove an $\epsilon$-Nash equilibrium property for the obtained decentralized strategies.
\end{abstract}

\maketitle

\section{Introduction}

\subsection{GMFG and GXMFG theory}

Graphon theory \cite{L12} provides a  natural framework for modeling the interaction of agents distributed over large networks when the nodes can be indexed by the reals in $[0,1]$ furnished only with the topology of the node set $[0,1]$. Graphons are    symmetric measurable functions on  $[0,1]^{2}$, taking values in $[0,1]$ (interpreted as generalized adjacency matrices). The graphon framework has been widely used in the analysis of dynamics, control and games on large networks  (see e.g. \cite{BCW20,CH21,
CCG22,ddP25,GC18, KX22,M14, PO19}).
 Graphon Mean Field Games (GMFGs) \cite{CH18,CH21} are formulated so as to describe the limiting behavior of large populations of dynamical agents interacting over large networks, the key property of interest being the existence of Nash equilibria.

The GMFG framework is formulated in terms of graphons  and hence is restricted to asymptotically dense networks. However, many networks of interest are not  dense.
Large sparse graphs occur, for instance, in  biological neural networks \cite{DM24} and electrical power grids \cite{PA13}.
When the sequence of graphs is sparse but the node degrees still have unbounded growth,  graphon theory is extended by an $L^p$ theory for sparse graph limits \cite{BCCZ19}.
The reader is also referred to sparse limit graph dynamical models in     \cite{LS21} for related game models on sparse networks, and  general convergence analyses of interacting particles over graphs, including sparse graphs,  in  \cite{KM17,LRW23}.
 The graphexon MFG (GXMFG) construction generalizes the GMFG framework to the graphexon setting (\cite{C22}) and consequently applies to networks which are asymptotically a combination of sparse and  dense graph limits. Unlike the graphon case,  non-trivial subsequential  graphexon limits are defined for all graph sequences and hence
in particular for sparse graph sequences.
The solution of the GXMFG \cite{CH24a} assigns coupling weights to other nodes using a graphexon measure instead of a graphon function in GMFG equations \cite{CHH22,CH21}.

\subsection{MFG and Laplexion dynamics}
Immediate neighbor based interactions have been addressed in the game theoretic  literature \cite{B93,I06,Y98} and arise in markets of retailing services or goods \cite{B93}. Such interaction structures are also adopted in deriving continuum dynamics as motivated by bacterial colonies and flocking \cite{FKS24}, and interacting particle models in mathematical physics \cite{BBP20,GPV88}.

To describe interactions of immediate neighboring clusters,
our previous work considers graphexon limits of sparse networks
and proposes first order interaction models \cite{CH24a} and second order interactions models \cite{CH24b} containing a mean field derivative term in the nodal parameter. This interaction structure accounts for the influence of  neighboring subpopulations through the spatial variations of their behaviors,  while the standard graphon network limit  framework has no defined metric topology on the set of vertices in which to define the notion of locality.
The more recent work \cite{CH25} starts with
a dynamical model via graph Laplacian operators and, taking their sparse network limit, provides a mathematical derivation of second order interaction terms in \cite{CH24b}, which, due to their interpretation via limit Laplacians, are called  Laplexions. Graph Laplacians have been widely used in other areas such as multi-agent coordination control \cite{OFM07} and machine learning \cite{CS20}.
For developing a limit theory, the modeling in this paper starts with finite populations on stylized finite networks as either one dimensional or two dimensional lattices \cite{BBP20,B93, FKS24,GPV88}.
Our network models are cases of extreme sparsity, and are not covered by the
sparse limit theory as in \cite{BCCZ19,KM17}.
In a linear-quadratic (LQ) setting, this paper addresses the relationship between the finite population model and  the sparse network limit model while some preliminary analysis has appeared in  \cite{CH25cdc} without  detailed proofs. Our analysis  depends on techniques  different from those employed in MFGs and GMFGs, primarily due to the second order (diffusive) interactions.
Our techniques will combine mean field approximations
for intra-cluster coupling and Laplexion approximations of the seconder order interactions for inter-cluster coupling. The Laplexion is closely related to hydrodynamic limit under scaling of diffusion type in asymptotic analysis of interacting particle systems \cite{BBP20,GPV88}.

The paper is organized as follows.
Section \ref{sec:sp} introduces the dynamical modeling involving a graph Laplacian operator to capture the influence of subpopulations  at neighbouring nodes, which leads to the so-called  Laplexion mean field game
(LMFG) in the large network limit. Section \ref{sec:lq} derives the solution equation system of the LQ LMFG and the associated decentralized strategies. Section \ref{sec:ex} presents an existence and uniqueness analysis of the partial differential equation (PDE) system of the LMFG. Section \ref{sec:ne} gives the a priori cost bound for the decentralized strategies.
Section \ref{sec:so} analyzes the Laplexion limit of the second order interaction terms.
Section \ref{sec:ne} establishes an asymptotic Nash equilibrium theorem.

{\bf Notation:} We denote by $N^\iota$ ($\iota$=1 for the ring, $\iota=2$ for the torus)  the number of nodes in the network, by $N_1$  the subpopulation size at each node,
by $C$ (or $C_1, \cdots$)  generic constants, which may change from place to place but do not depend on $(N,N_1)$, and by $I_k$ the $k\times k$ identity matrix.

\section{Dynamic games with subpopulations over sparse networks and the Laplexion limit}
\label{sec:sp}

We start with some preliminaries on graph limits and graph Laplacian operators.

\subsection{Large graphs and their limits}

Let $(\mathcal{M}, {\mathbf g})$ be a compact, connected Riemannian manifold of dimension $p$, equipped with its geodesic distance $d_{\mathbf g}$ and Riemannian volume measure $dV_{\mathbf g}$.
Let ${\mathcal P}({\mathcal M})$ be all probability measures defined on $({\mathcal M}, {\mathcal B}({\mathcal M}))$ with Borel $\sigma$-field ${\mathcal B}({\mathcal M})$.
\begin{definition}[Embedded graph sequence]\label{def:embedded_graph_sequence_simple}
 $G^N=(\mathcal{V}^N,\mathcal{E}^N)$ shall denote a non-empty finite graph without self-loops, with node set denoted $\mathcal{V}^N  \coloneqq \{1,\cdots,N\}$, and undirected edge set denoted $\mathcal{E}^N \subseteq \mathcal{V}^N \times \mathcal{V}^N$. An edge between two distinct nodes $i$ and $j$ is denoted by the unordered pair $(i;j)$.
The sequence $\{G^N,\varphi^N\}_{N=1}^\infty$ is an \emph{embedded graph sequence} if every $G^N$ is embedded into $\mathcal{M}$ via a bijective node mapping $\varphi^N: \mathcal{V}^N \rightarrow \mathcal{M}$, assigning a location $\alpha_i^N = \varphi^N(i) \in \mathcal{M}$ to each node $i \in \mathcal{V}^N$.
\end{definition}

\begin{definition}[Empirical vertex and edge measures]\label{def:empirical_measures}
Let $G_N=(\mathcal{V}_N,\mathcal{E}_N)$ be a finite embedded graph on a Riemannian manifold ${\mathcal M}$ with an embedding mapping $\varphi^N:\mathcal{V}_N\to {\mathcal M}$, with $\alpha_i^N=\varphi^N(i)$.
\begin{enumerate}
    \item The \emph{empirical vertex measure} $V_N \in {\mathcal P}({\mathcal M})$ is defined by
    \[
        V_N \; \coloneqq  \; \frac{1}{|\mathcal{V}_N|}\sum_{i\in\mathcal{V}_N}\delta_{\alpha_i^N},
    \]
    where $\delta_{\alpha_i^N} \in \mathcal{P}({\mathcal M})$ denotes the Dirac probability measure on $({\mathcal M},\mathcal{B}({\mathcal M}))$ concentrated on $\alpha_i^N$.

\item The \emph{empirical graphexon measure} $W_N \in {\mathcal P}({\mathcal M}\times{\mathcal M})$ is defined by
\[
  W_N \; \coloneqq \; \frac{1}{2|\mathcal{E}^N|}\sum_{(i;j)\in \mathcal{E}^N}\big(\delta_{(\alpha_i^N,\alpha_j^N)}+
  \delta_{(\alpha_j^N,\alpha_i^N)}\big),
\]
where $\delta_{(\alpha_i^N, \alpha_j^N)}\in {\mathcal P}({\mathcal M}\times{\mathcal M})$ denotes the Dirac probability measure on $({\mathcal M}\times {\mathcal M},\mathcal{B}({\mathcal M}\times {\mathcal M}))$ that concentrates on the ordered pair $(\alpha_i^N,\alpha_j^N)$ whenever $(i;j)\in\mathcal E^N$.
\end{enumerate}
\end{definition}
\begin{remark}
     The measure $W_N$ assigns zero measure to the diagonal in ${\mathcal M} \times {\mathcal M}$ since the graphs under consideration have no self-loops.
\end{remark}

Since $\mathcal{M}$ is a compact metric space, the space of probability measures ${\mathcal P}(\mathcal{M})$ is tight. Consequently, by Prokhorov's theorem, the space of probability measures ${\mathcal P}(\mathcal{M})$, endowed with the topology of weak convergence, is also compact. Consequently, the sequence of empirical vertex measures $\{V_N\}_{N \in \mathbb{N}}$ is relatively compact which guarantees the existence of convergent subsequences.

\begin{definition}[\cite{C22,CH24a,CH24b}]\label{def:vertexon_graphexon_limits}
Let $(V_N,W_N)$ be the empirical vertex and edge measures from Definition~\ref{def:empirical_measures}.
We equip ${\mathcal P}({\mathcal M})$ with the topology of weak convergence of probability measures, and let
$\mathcal{P}_s({\mathcal M}\times{\mathcal M})$ denote the space of finite \emph{symmetric} Borel measures on ${\mathcal M}\times{\mathcal M}$.
\begin{itemize}
    \item A \emph{vertexon} is any weak limit $V_\infty\in {\mathcal P}({\mathcal M})$ of a subsequence $\{V_{N_k}\}_{k=1}^\infty$.

    \item Given such a vertexon $V_\infty$ realized by $\{N_k\}_{k=1}^\infty$, an associated \emph{graphexon} is any weak limit
  $W_\infty\in\mathcal{P}_s({\mathcal M}\times{\mathcal M})$ of a further   subsequence in $\{W_{N_k}\}_{k=1}^\infty$.
\end{itemize}
\end{definition}

The limit $V_\infty$ depends on the chosen subsequence
$\{N_k\}_{k=1}^\infty$.
 Hence, both $V_\infty$ and the associated $W_\infty$ need not be unique.

\subsection{The graph Laplacian operator}
\label{sec:sub:graph}

Consider an undirected graph $G=({\cal V}, {\cal E})$ without self-loops, where ${\cal V}=\{\alpha_1, \cdots, \alpha_N\}$ is the set of (distinct) nodes and ${\cal E}$ is the set of edges. An edge between $\alpha_i$ and $\alpha_j$, $i\ne j$, is denoted by an unordered pair $(\alpha_i;\alpha_j)$;
thus both $(\alpha_i;\alpha_j)$ and $(\alpha_j;\alpha_i)$ mean the same edge. For convenience in the subsequent analysis, we label the node set in this form instead of using a set of integers.
The finite graph Laplacian $(g_{\alpha_i \alpha_j})_{1\le i,j,\le N}$ is a symmetric matrix, where for $i\ne j$, $g_{\alpha_i \alpha_j}>0$ if and only if $(\alpha_i;\alpha_j)$ is an edge (in which case we denote $\alpha_i\sim \alpha_j$), and $g_{\alpha_i \alpha_j}=0$ otherwise; $g_{\alpha_i \alpha_i}= -\sum_{j \ne i} g_{\alpha_i \alpha_j}$.

 For a function $\psi$ from ${\cal V}$ to ${\mathbb R}^n$, denote the finite graph Laplacian operator by
$$
({\cal L}_N \psi )({\alpha_i}) := \sum_{\alpha_j: \alpha_j\sim \alpha_i }g_{\alpha_i \alpha_j}(\psi({\alpha_j})-\psi({\alpha_i})).
$$
So ${\cal L}_N \psi$ is still a function from  ${\cal V}$
to ${\mathbb R}^n$.

\subsection{Agent dynamics and costs within a  finite population}

For the graph   $G$   in Section \ref{sec:sub:graph}, suppose that
each of the $N$ nodes is  occupied by a cluster (or called a subpopulation) of agents.
At  node $\alpha_k$, $1\le k\le N$,  denote the agent states by $x_{\alpha_k}^i(t)$, $1\le i\le N_1$. For notational simplicity, we assume equal subpopulation sizes $N_1$ although this restriction can be relaxed.
Define the local empirical state distribution of the subpopulation at node $\alpha_k$:
 \begin{align}
 \mu_{\alpha_k}^{(N_1)}(t)= \frac{1}{N_1}\sum_{j=1}^{N_1}\delta_{x_{\alpha_k}^j(t)},
 \quad 1\le k\le N, \notag
 \end{align}
 where $\delta_x$ is the Dirac measure at the point $x$.

While we are primarily interested in developing our methodology in an LQ setting, our exposition starts with dynamics and costs in a general form to make the notation more compact.
Let $(\Omega, {\mathcal F}, \{\mathcal F_t\}_{t\ge 0}, {\mathbb P})$ be the underlying filtered probability space satisfying the usual condition.
For agent dynamics we introduce the stochastic differential equation (SDE):
\begin{align}
d x_{\alpha_k}^i(t)&= f(x_{\alpha_k}^i(t),u_{\alpha_k}^i(t),\mu_{\alpha_k}^{(N_1)}(t)) dt\notag \\
&\quad +\sum_{m: m\ne k}g_{\alpha_k \alpha_m}\Big[ \int_{{\mathbb R}^n} \phi (y)\mu_{\alpha_m}^{(N_1)}(t,dy) - \int_{{\mathbb R}^n} \phi (y)\mu_{\alpha_k}^{(N_1)}(t,dy)\Big]dt+\sigma dw_{\alpha_k}^i, \label{dxiakN}
\end{align}
where the state $x_{\alpha_k}^i(t)$, control $u^i_{\alpha_k}(t)$ and standard Brownian motion $w_{\alpha_k}^i(t)$ have values in ${\mathbb R}^n$, ${\mathbb R}^{n_1}$, and ${\mathbb R}^{n_2}$, respectively; $\phi=(\phi_1, \cdots \phi_n)^T: {\mathbb R}^n \to {\mathbb R}^n$  is a given coupling function; and $\sigma$ is a constant coefficient matrix.    The ${\mathcal F}_0$-measurable initial states  $x_{\alpha_k}^i(0)$, $1\le i\le N_1, 1\le k\le N$, are independent. The ${\mathcal F}_t$-adapted  standard Brownian motions  are independent and also independent of the initial states.  The above graph Laplacian based inter-cluster coupling is different from the dynamical coupling in \cite{CH21} since in the present model we do not assume a dense graph sequence.

The individual cost is given by
\begin{align}
J_{\alpha_k}^i= {\mathbb E}\int_0^T L(x_{\alpha_k}^i(t),u_{\alpha_k}^i(t),  \mu_{\alpha_k}^{(N_1)} (t))dt + \EE h(x_{\alpha_k}^i(T),\mu_{\alpha_k}^{(N_1)} (T) ). \label{cjikN}
\end{align}

\subsection{Limit dynamics of infinite sub-populations on a finite graph}

As an approximation of \eqref{dxiakN}--\eqref{cjikN} with
$N_1\to \infty$,
we next introduce  a mid-level model where the  $N$ nodes of the graph $G$ are each occupied by  an infinite subpopulation.  At time  $t$, the agents of the cluster residing at node $\beta\in {\cal V}=\{\alpha_1, \cdots,\alpha_N\}$  generate the local mean field  $\mu_\beta(t)$, as the infinite population limit of $\mu_\beta^{(N_1)}(t)$ when $N_1\to \infty$.

Given the coupling function $\phi$, define  $\bar \phi_i(t,\beta)=\int_y \phi_i (y) \mu_\beta (t,dy)\in {\mathbb R}$, $1\le i\le n$, and
$\bar \phi=(\bar \phi_1,\cdots, \bar\phi_n)^T$  for $\beta\in {\cal V}$. Let $x_\alpha(t)$ denote the state of a  representative agent  at node $\alpha\in {\mathcal V}$.
 We consider the following  dynamics
\begin{align}\label{sdeL}
dx_\alpha(t)=f(x_\alpha(t), u_\alpha(t), \mu_\alpha(t))dt  +
({\cal L}_N \bar \phi( t,\cdot))(\alpha) dt + \sigma dw_\alpha(t) ,
\end{align}
where  ${\cal L}_N$ is the finite graph Laplacian operator.

We further introduce the cost of the representative agent
\begin{align}
J_\alpha= \mathbb{E}\int_0^T
L(x_\alpha(t), u_\alpha (t) , \mu_\alpha(t)) dt+\EE h(x_\alpha(T), \mu_\alpha(T)). \notag
\end{align}

\subsection{The solution via infinite populations and the large network limit}

\label{sec:rtop}

\begin{figure}[t]
\begin{center}
\psfig{file=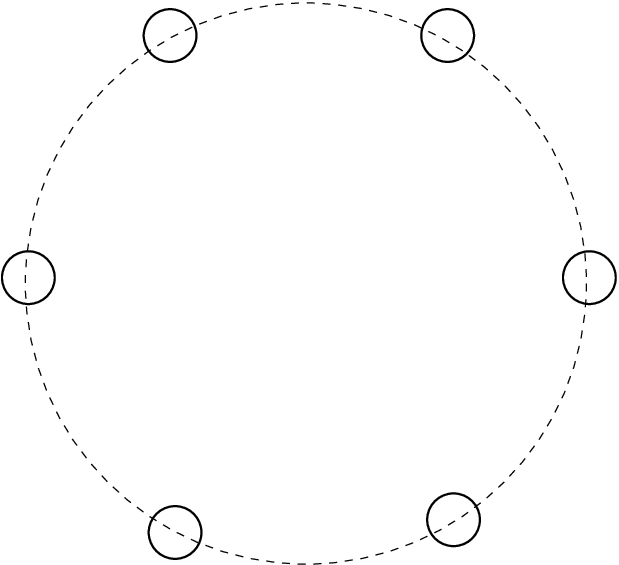, width=3.4cm, height=3cm}
\end{center}
\caption{The ring network with two immediate neighbors  }
\label{fig:ring}
\end{figure}

For asymptotic analysis with a suitable neighbor structure,  we take the ring network ${\mathbb T}=[0,1)$ (where 1 is identified with 0) and the torus network ${\mathbb T}^2$ as prototype network limit models.
They are intended to approximate large one and two dimensional networks while their periodic structure avoids the nuisance of boundary effects.
\begin{itemize}

\item[(I)] {\bf Discrete ring network} ${\mathbb T}_{1/N}$.
Consider the partition of $[0,1)$ into $N$ subintervals
${\cal I}_k=[(k-1)/N, k/N)$, $k=1, \cdots, N$.
Let the midpoint $\alpha_k=\frac{k-1}{N} +\frac{1}{2N} $ denote the nodal location.    Denote the nodal set ${\cal V}_{1/N}=\{\alpha_k, 1\le k\le N  \} $ forming a ring, as a graph, denoted by ${\mathbb T}_{1/N}$ (see Fig. \ref{fig:ring}) so that $\alpha_1\sim \alpha_N$. We define the graph Laplacian as an $N\times N$ matrix: for $j\ne i$,
$g_{\alpha_i \alpha_j}=
N^2$ if $ \alpha_j\sim \alpha_i,
$
and $g_{\alpha_i \alpha_j} =0$ otherwise;
and $g_{\alpha_i\alpha_i}= -2N^2 $.
\item[(II)] {\bf Discrete torus network} ${\mathbb T}^2_{1/N}$.
The discrete torus is a graph consisting of
$N^2$ nodes of the form $\alpha_{ij}\coloneqq (\alpha_i, \alpha_j)$ (an ordered pair), where $\alpha_i, \alpha_j\in {\cal V}_{1/N}$, so that each node $\alpha_{ij}$ has $4$ neighbors. The nodal set is
${\cal V}^2_{1/N} $. The graph Laplacian is defined by the rule: for $\alpha_{ij}\ne\alpha_{i'j'} $,  $g_{\alpha_{ij}\alpha_{i'j'}}=N^2$ if $\alpha_{ij}\sim \alpha_{i'j'} $, and $g_{\alpha_{ij}\alpha_{i'j'}}=0$ otherwise; and $g_{\alpha_{ij}\alpha_{ij}}=-4N^2$ for all $i,j$.

\end{itemize}

In the mathematical physics literature,
a scaling by $N^2$ for immediate neighbors is customarily called
the scaling of diffusion type \cite{GPV88}, which is used to derive the so-called hydrodynamic limits.
We will analyze  limits of the finite graph Laplacian operators in \eqref{sdeL}.

Recall that $\bar \phi_i(t, \alpha)=\int \phi_i(y) \mu_\alpha(t,dy) $ for $\alpha\in {\cal V}$. To analyze the large network limit, we take the  continuum vertex set ${\mathcal V}_\infty$ as  ${\mathbb T}^\iota$, $\iota=1,2$, and let $\mu_\alpha(t) $ be defined for all $\alpha\in {\cal V}_\infty$, so that $\bar \phi_i(t, \alpha)$ is now defined for $\alpha\in {\cal V}_\infty=[0,1)^\iota$. We further identify $\bar \phi_i$ as a function defined on $[0,T]\times \mathbb{R}^\iota$ (still denoted as $\bar \phi_i$), with period 1 for each component of ${\mathbb{R}^\iota}$.
A function $f(t,z)$ on $[0,T]\times \mathbb{R}^k$, $k\ge 1$, is said to have period 1 for each space variable if
$f(t, z_1+i_1, \cdots z_k+i_k)= f(t,z_1, \cdots, z_k)$ for all $z_1, \cdots, z_k\in {\mathbb R}$ and all $i_1, \cdots, i_k\in {\mathbb Z}$.
If $\bar \phi_i$ has up to second order continuous partial derivatives in the variable $\alpha\in {\mathbb R}^\iota$, by elementary calculation,  we derive
$$
\lim_{N\to \infty}\sup_{\alpha\in {\cal V}^\iota_{1/N}}|{\cal L}_N \bar \phi(t,\cdot)(\alpha) -\Delta\bar \phi(t,\alpha)|=0.
$$
where
\begin{align*}
\Delta\bar \phi(t,\alpha)&\coloneqq (\Delta \bar\phi_1(t,\alpha),\cdots,\Delta \bar\phi_n(t,\alpha) )^T, \\
  \Delta \bar \phi_i(t,\alpha)&=
\begin{cases}
 \partial_\alpha^2\bar \phi_i(t,\alpha) & \mbox{if} \  \alpha\in {\mathbb T},   \\
\sum_{k=1}^2\partial_{\alpha^{(k)}}^2  \bar \phi_i (t,\alpha) & \mbox{if} \ \alpha=(\alpha^{(1)}, \alpha^{(2)})\in {\mathbb T}^2.
\end{cases}
\end{align*}

The following limit
model for \eqref{sdeL} takes a formal limit of the graph Laplacian operator and  is called a
{\it {Laplexion Mean Field  Dynamical System}}:
\begin{align}\label{sdeLJ}
dx_\alpha(t)&=f(x_\alpha(t), u_\alpha(t), \mu_\alpha(t))dt  +\bar \Phi(t,\alpha) dt   + \sigma dw_\alpha(t),
\end{align}
which specifies the dynamics of a representative agent $ {\mathcal A}_\alpha$ at node $\alpha \in {\cal V}_\infty$ and where we suppose
$\bar \Phi(t,\alpha) \coloneqq \Delta\bar \phi(t,\alpha)$ is well-defined.
 The cost function is given by
\begin{align}
J_\alpha&= \mathbb{E}\int_0^T
L(x_\alpha(t), u_\alpha (t) , \mu_\alpha(t)) dt +\EE h(x_\alpha (T), \mu_\alpha(T)).
\end{align}

The above Laplexion framework is  different from the GMFG case which depends on a graphon function for coupling.
 For the ring case of ${\mathbb T}_{1/N}$, the graphexon $W_\infty$, as  the counterpart of a graphon, assigns a uniform probability distribution along the diagonal line of the unit square $[0,1)^2$, and as a result, the use of $W_\infty$ alone cannot indicate interconnections between clusters. By contrast,
the incorporation of the Laplacian operator enables us to capture the immediate neighbor interactions in the  limit model \eqref{sdeLJ}.

\subsection{The Laplexion mean field game (LMFG) equations }

\label{gmfg}

Denote   $\Sigma_w=\sigma \sigma^T $.   Based on locally infinite populations and the network limit, the formal LMFG solution system consists of the pair of
  Hamilton--Jacobi--Bellman (HJB) equation and
 Fokker--Planck--Kolmogorov (FPK) equation:
\begin{align}
-\partial_t V_\alpha (t,x)& =\inf_u
\Big\{\partial_x V_\alpha(t, x) [ f ( x, u, \mu_\alpha(t)) +
\bar \Phi(t,\alpha )]  + L(x, u, \mu_\alpha(t)) \Big\}  \label{hjbPh} \\
&\qquad\quad +\frac{1}{2} \mbox{Tr} (\Sigma_w \partial_x^2 V_\alpha (t, x)),\notag  \\
 V_{\alpha  }  (T,x)& =h(x,\mu_\alpha(T)),    \quad(t,x)\in [0,T]\times\mathbb{R}^n, \quad  \alpha \in {\mathcal V}_\infty, \notag\\
 \label{FPK}
{\partial_t p_\alpha(t,x)} & =-\mbox{div} \{[ f( x,    u^0, \mu_\alpha(t)) +\bar \Phi(\alpha,t) ] p_\alpha(t,x) \}   + \frac{1}{2}
\mbox{Tr} (\Sigma_w  {\partial^2_x p_\alpha(t,x)}),
\end{align}
where $u^0\coloneqq  u^0(t, x| {\mu}_{G } ) $ is the best response control law determined by equation \eqref{hjbPh},
$\bar \Phi(t,\alpha )=\Delta \int_{{\mathbb R}^n} \phi(y) p_\alpha (t,y)dy $ with coupling function $\phi$, and $\mu_G=(\mu_\beta)_{0\le\beta < 1}$.
We use $p_\alpha(t, \cdot)$ to denote the probability density function of  $\mu_\alpha(t)$ (as the local mean field at node $\alpha$) if it exists. The term $\partial_x^2 V_\alpha$ stands for the Hessian matrix, and $\mbox{div}(\cdot)$ for divergence.

As in standard mean field games,
the HJB equation \eqref{hjbPh} is derived by solving an optimal control problem of agent ${\mathcal A}_\alpha$  taking $\mu_\alpha(\cdot)$  as fixed.
The solution of the HJB--FPK equation system then consists of
the family of function pairs
$(V_\alpha, p_\alpha )$ defined on $[0,T]\times \mathbb{R}^n$, which are indexed by $\alpha \in {\mathcal V}_\infty$.

\section{The Linear Quadratic  (LQ) Model}

\label{sec:lq}

In all subsequent sections,  we will study the model \eqref{dxiakN}--\eqref{cjikN}
in an LQ setting.
 The cases where the finite network   is respectively either a discrete ring
 ${\mathbb T}_{1/N}$ or a discrete torus ${\mathbb T}^2_{1/N}  $ will be considered: first,   the ring has nodes $\{\alpha_1^N, \cdots,\alpha_{N^*}^N  \}$, where $ N^*= N$, while, second, the torus network is such that $ N^*= N^2 $ .   Let $N_1$ be the equal subpopulation size at all nodes.   Denote the subpopulation average state
$  x_{\alpha_k^N}^{(N_1)} = \frac{1}{N_1} \sum_{j=1}^{N_1} x^j_{\alpha_k^N}.$
The agent ${\cal A}^i_{\alpha_k^N}$, $1\le i\le N_1$, at node $\alpha_k^N $ has the following dynamics:
\begin{align}
dx^i_{\alpha_k^N}(t) =\ &[ Ax^i_{\alpha_k^N}(t)
+ Bu^i_{\alpha_k^N}(t)
+ D  x_{\alpha_k^N}^{(N_1)}(t)]dt \label{dxat3} \\
&+
\kappa N^2 \sum_{k':\alpha_{k'}^N\sim \alpha_k^N } (x_{\alpha_{k'}^N}^{(N_1)}
-x_{\alpha_k^N}^{(N_1)})(t)dt +\sigma dw_{\alpha_k^N}^i(t),\notag
\end{align}
with initial state $x^i_{\alpha_k^N}(0)  $,
where $N^2\sum_{k': \alpha_{k'}^N\sim \alpha_k^N} (\cdot) $ is called the second order interaction term, involving  the subpopulations
at  immediate neighboring nodes (see Fig. \ref{fig:ring} for the ring network  case). Here $x^i_{\alpha_k^N}(t)$, $ u^i_{\alpha_k^N}(t)  $, and $w^i_{\alpha_k^N}(t)  $ have dimensions $n$, $n_1$, and $n_2$, respectively, and $\kappa>0$ is the interaction coefficient.
 For  a matrix $M\ge 0$, denote $\|z\|_M^2= z^TMz$.
The cost of agent ${\cal A}^i_{\alpha_k^N}$  is
\begin{align}
 &J_{\alpha_k^N}^i(u^i_{\alpha_k^N}, u^{-i}_{\alpha_k^N},  {\bf u}_{-\alpha_{k}^N}  ) \label{jakiu} \\
 & = \EE \int_0^T\Big[ \|x_{\alpha_k^N}^i(t)-\Gamma_1 x^{(N_1)}_{\alpha_k^N}(t)-\Gamma_2\eta\|_Q^2   + \|u_{\alpha_k^N}^i(t)\|_R^2 \Big]dt \notag \\
&\qquad\quad + \EE  \|x_{\alpha_k^N}^i(T)-\Gamma_{1f} x^{(N_1)}_{\alpha_k^N}(T)-\Gamma_{2f}\eta\|_{Q_f}^2,
\notag
\end{align}
for $R>0$, $Q\ge 0$, $Q_f\ge 0$, and $\eta\in \mathbb{R}^n$,
 where $u_{\alpha_k^N}^{-i}$ stands for the strategies of all agents at node $\alpha_k^N$ except ${\cal A}_{\alpha_k^N}^i$, ${\bf u}_{-\alpha_{k}^N}$  for strategies of all agents not at node $\alpha_k^N$.
 For notational simplicity, we consider constant parameters $(A,B, \cdots, \eta, Q_f)$. The extension to time-dependent parameters will be straightforward.

Within the large network limit, we take ${\mathbb T}^\iota$ ($\iota=1, 2$) as the node set, and $x_\alpha(t)\in \mathbb{R}^n$ as the state of the representative agent  ${\cal A}_\alpha$ situated at node $\alpha\in
{\mathbb T}^\iota $.
At node $\alpha$, the state average of the subpopulation is replaced by $m(t,\alpha)=\int_{{\mathbb R}^n} y \mu_\alpha (t,dy).
$
We take  $\phi(y)=\kappa y$, $\kappa>0$,
so that $\bar \Phi(t,\alpha)$ in \eqref{sdeLJ}
becomes $ \kappa \Delta m(t,\alpha) $.
Suppose all initial states have  mean $m_\alpha(0)$
and finite second moment at $t=0$.

For the best response control problem of agent  ${\cal A}_\alpha$, we consider the dynamics \begin{align}
dx_\alpha(t)    = &
(A x_\alpha(t) +Bu_\alpha(t) +D m(t,\alpha)
+\kappa \Delta m(t,\alpha) ) dt+   \sigma dw_\alpha (t).
 \label{lqdyn}
\end{align}
The cost of agent ${\cal A}_\alpha$ is given by
\begin{align}
J_\alpha &= {\mathbb E} \int_0^T\Big( ||x_\alpha(t)-\Gamma_1 m(t,\alpha)-\Gamma_2 \eta ||_{ Q}^2+ u_\alpha^T(t) R u_\alpha(t)\Big)dt \notag  \\
&\quad + {\mathbb E} \|x_\alpha(T)-\Gamma_{1f}m(T,\alpha)-\Gamma_{2f}\eta\|^2_{Q_f}.\label{lqcos}
\end{align}

Next, we introduce the Riccati ordinary differential equation (ODE):
\begin{align*}
0=\dot{P}+  PA +A^T P-PBR^{-1} B^T P +Q,\qquad P(T)=Q_f,
\end{align*}
 which has  a unique solution $P$ on $[0,T]$.
The best response control law is given by
\begin{align}
\hat u_\alpha(t)=  -R^{-1} B^T
[P(t) x_\alpha (t) +S(t,\alpha) ],
\label{bru}
\end{align}
where $S(t,\alpha)$ is given by an ODE when $\alpha$ and the function $m$ are fixed.
Subsequently, $m$ is determined by the standard  consistency condition of MFGs \cite{HZ20}: the ODE of $m (t,\alpha)$  is obtained by averaging the closed-loop equations
for \eqref{lqdyn} of all agents residing at node $\alpha$.
We summarize the LMFG equation system:
\begin{align}
\partial_t{S}(t,\alpha)&= -(A^T  -P(t)BR^{-1}B^T)S(t,\alpha)  \label{pde4S}  \\
&  \quad -\kappa P(t)  \Delta m(t,\alpha) +(Q\Gamma_1-P(t)D)  m(t,\alpha) + Q\Gamma_2 \eta, \notag \\
\partial_t{m (t,\alpha)} &=\kappa\Delta m(t,\alpha) + (A-BR^{-1}B^T P(t)+D )m(t,\alpha )\label{pde4m}  \\
&\quad  - B R^{-1} B^T S(t,\alpha), \qquad  0<t<T, \quad \kappa>0,    \notag
\end{align}
where $\alpha\in {\mathbb T}^\iota$,   $m(0,\alpha)$ is given, and
$
S(T,\alpha)= -Q_f (\Gamma_{1f}m(T,\alpha)+\Gamma_{2f}\eta   ) .
$
The equation system  is  a generalization of the
 forward-backward ODE system (see e.g. \cite[Sec. IV]{HZ20}) in MFGs with symmetric players.
The condition $\kappa>0$ is important in existence analysis since \eqref{pde4m} becomes a parabolic equation (with space variable $\alpha$) subject to an initial condition.

\section{Existence and Uniqueness Analysis}
\label{sec:ex}

To apply results in the literature to the existence analysis
for \eqref{pde4S}--\eqref{pde4m}, we will identify $S$ and $m$ as functions defined on $[0,T]\times \mathbb{R}^\iota    $ ($\iota=1$ or 2) with period 1 in each space variable within $ \alpha$.

For a  general analysis of classical solutions,  we introduce   certain H\"older spaces \cite{F08,LUS68}. Fix $\delta \in (0,1)$ as the H\"older exponent.
On  ${\mathbf Q}= [0,T] \times {\mathbb R}^k$ with $k\ge 1$, we define the parabolic distance between two points
$(t,x)$ and $(s,y)$ as
$$
d_p((t,x), (s,y))=  |t-s|^{1/2} +\max_{i\le k}|x_i-y_i|.
$$
  Define the H\"older semi-norm for a function $v(t,x)$ with domain ${\mathbf Q}$:
$$[v]_{\delta/2, \delta;{\mathbf  Q}}=\sup_{z_i\in {\mathbf Q}, z_1\ne z_2} \frac{|v(z_1)-v(z_2)|}{d_p(z_1, z_2)^\delta}.
$$
We allow $v$ to be a vector function.
Denote
$|v|_{0;{\bf Q}}= \sup_{z\in {\bf Q}}|v(z)|$.
 We further define the H\"older norms
\begin{align*}
|v|_{\delta/2, \delta;{\mathbf  Q}}= \ &|v|_{0;{\bf Q}} +[v]_{\delta/2, \delta;{\mathbf Q}},\\
|v|_{1+\delta/2, 2+\delta; {\mathbf Q}}=\ & |v|_{0;{\bf Q}}+|\partial_tv|_{0;{\bf Q}} +\sum_i|\partial_{x_i} v|_{0;{\bf Q}}\\
& + \sum_{i,j}|\partial_{x_ix_j} v|_{0;{\bf Q}} + [v]_{1+\delta/2, 2+\delta; {\mathbf Q}},
\end{align*}
 where
$[v]_{1+\delta/2, 2+\delta; {\mathbf Q}}= [\partial_t v]_{\delta/2, \delta; {\mathbf Q}}+\sum_{i,j} [\partial_{x_ix_j} v]_{\delta/2, \delta; {\mathbf Q}}$.
We use $C^{\delta/2, \delta}({\mathbf Q})  $  (resp., $C^{1+\delta/2, 2+\delta}({\mathbf Q})$) to denote the space consisting of functions $v$ with $|v|_{\delta/2, \delta; {\mathbf Q}}<\infty$
(resp., $|v|_{1+\delta/2, 2+\delta; {\mathbf Q}}<\infty$).
The H\"older space $C^{l+\delta}({\mathbb R}^k)  $ for $l=0,1, 2,\cdots$,  is similarly defined without involving time variable $t$.

We give some analytical preparation on linear parabolic equations.
Let $q(t,z)$ be an $\mathbb{R}^n$-valued function defined on ${\mathbf Q}$.
Consider the following parabolic equation system \cite[Theorem 5.1, Theorem 10.2]{LUS68}:
\begin{align}
\begin{cases}
\partial_t q(t,z) =\kappa \Delta q(t,z) +F(t,z)
q(t,z )+g(t,z), \\
 q(0,z )=\varphi(z),
\end{cases} \label{auxpde}
\end{align}
where $\kappa>0$,
$\varphi \in C^{2+\delta}({\mathbb R}^k)$,
$F, g\in C^{\delta/2, \delta}({\bf Q})$,
and
$|F|_{\delta/2, \delta; {\mathbf Q}} \le C_F. \notag$
Here $F$ and $g$ are $\mathbb{R}^{n\times n}$-  and $\mathbb{R}^n$-valued, respectively.

The following lemma is a special case of Theorem 10.2 in \cite{LUS68}.
For notational brevity, we will omit the underlying domain for the norms.

\begin{lemma}\cite{LUS68}  \label{lemma:aux}
Suppose ${\mathbf Q}=[0,T]\times {\mathbb R}^k$.
Then there exists a unique  solution  $ q\in  C^{1+\delta/2,2+\delta}({\mathbf Q})$ to \eqref{auxpde} and
\begin{align*}
|q|_{1+\delta/2, 2+\delta} \le
C_0 (|g|_{\delta/2,\delta} + |\varphi|_{ 2+\delta} ),
\end{align*}
where the constant $C_0$ only depends on
$(\kappa, C_F, T, \delta,k,n)$.
\end{lemma}

\begin{remark}
$C_0$ in Lemma \ref{lemma:aux} can be selected such that it remains valid when $T$ is replaced by any $T'<T$.
\end{remark}

For existence analysis,
we introduce the following function classes (with $\iota= 1,2$):

\begin{itemize}

\item[(i)] ${\cal C}_{prd}^\delta({\mathbb R}^\iota)$ consists of functions $f\in C^\delta({\mathbb R}^\iota)$, having period 1 for each space variable.

\item[(ii)] ${\cal C}_{prd}^{2+\delta}({\mathbb R}^\iota)$ consists of functions $f\in C^{2+\delta}({\mathbb R}^\iota)$, having period 1 for each space variable (thus $\partial_{z_iz_j}f(z)\in C^\delta({\mathbb R}^\iota)$ for all $i,j$).
$C_{prd}^{4+\delta}({\mathbb R}^\iota)$ is similarly defined.

\item [(iii)] $ {\cal C}_{prd}^{\delta/2,\delta}( [0,T] \times   {\mathbb R}^\iota   )$  consists of functions $f\in C^{\delta/2,\delta}([0,T]\times{\mathbb R}^\iota)$, having  period 1 for each space variable.

\item [(iv)] $ {\cal C}_{prd}^{1+\delta/2,2+\delta}(  [0,T] \times  {\mathbb R}^\iota   )$  consists of functions $f$ in $C^{1+\delta/2,2+\delta}(  [0,T]\times  {\mathbb R}^\iota  )$, having  period 1 for each space variable.

\end{itemize}

For convenience of subsequent analysis, we introduce another H\"older semi-norm.
For two points
$x$ and $y$ in $\mathbb{T}^\iota =[0, 1)^\iota$,  define the distance
$d_{{\mathbb T}^\iota}(x,y)=\max_{i\le \iota} \min\{ 1-|x_i-y_i|,  |x_i-y_i|\}.$
Next, on the set ${\mathbf Q}^\circ= [0,T]  \times {\mathbb T}^\iota$, we define the parabolic distance between two points
$(t,x)$ and $(s,y)$ as
$$
d_p^\circ((t,x), (s,y))= |t-s|^{1/2} + d_{{\mathbb T}^\iota}(x,y).
$$
The H\"older semi-norm for a function $v(t,x)$ defined on ${\mathbf Q}^\circ$ is given by
$$
[v]_{\delta/2, \delta;{\mathbf  Q}^\circ}=\sup_{z_i\in {\mathbf Q}^\circ, z_1\ne z_2} \frac{|v(z_1)-v(z_2)|}{d_p^\circ(z_1, z_2)^\delta}.
$$
Denote
$|v|_{0;{\bf Q}^\circ}= \sup_{z\in {\bf Q}^\circ}|v(z)|
$ and further define the H\"older space $C^{\delta/2,\delta}({\mathbf Q}^\circ )  $.

\begin{lemma}\label{lem:qholder}
Suppose $v\in C_{prd}^{\delta/2, \delta}({\mathbf Q})$ with ${\mathbf Q}=[0,T]\times {\mathbb R}^\iota$.
Then we have
\begin{align}
[v]_{\delta/2, \delta;{\mathbf  Q}}= [v]_{\delta/2, \delta;{\mathbf  Q}^\circ} .
\end{align}
\end{lemma}
\begin{proof}
We start by showing $[v]_{\delta/2, \delta;{\mathbf  Q}}\le [v]_{\delta/2, \delta;{\mathbf  Q}^\circ}$.
 Take  two different points $(t,x),(t',x')\in {\mathbf  Q} $, with $x, x'\in {\mathbb R}^\iota $. Using periodicity,  we can find $x''\in {\mathbb R}^\iota$ such that  i) both $x''$ and $x'$ are identified with the same point in ${\mathbb T}^\iota =[0,1)^\iota$, ii)  $v(t',x')=v(t',x'')$, and moreover, iii) $|x_i-x''_i|\le \min \{ |x_i-x'_i|,1/2\}$ for each $i\le \iota$. Then we have
 \begin{align}
\frac{|v(t,x)-v(t',x')|}{d_p((t,x),(t',x'))^\delta }\le
\frac{|v(t,x)-v(t',x'')|}{d_p((t,x),(t',x''))^\delta } =
 \frac{|v(t,x)-v(t',x'')|}{d^\circ_p((t,x),(t',x''))^\delta }   \le [v]_{\delta/2, \delta;{\mathbf  Q}^\circ}.\notag
 \end{align}
 Hence  $[v]_{\delta/2, \delta;{\mathbf  Q}}\le [v]_{\delta/2, \delta;{\mathbf  Q}^\circ}$. We can similarly show $[v]_{\delta/2, \delta;{\mathbf  Q}^\circ}
 \le   [v]_{\delta/2, \delta;{\mathbf  Q}} $.
\end{proof}

Suppose $m(0,\cdot)\in C_{prd}^{2+\delta}({\mathbb R}^\iota)$ in \eqref{pde4m}.
We proceed to the existence analysis of \eqref{pde4S}--\eqref{pde4m}
by formulating the following fixed point problem with ${\mathbf Q}=[0,T]\times {\mathbb R}^\iota$.
Given any ${\mathbb R}^n$-valued function  $ S\in C_{prd}^{\delta/2, \delta}({\mathbf Q}) $  in \eqref{pde4m}, by Lemma
\ref{lemma:aux}, we obtain   a unique  solution
$ m$ and  define the following mapping $\Lambda_1$:
$$
 m =\Lambda_1( S),
$$
which is well defined and is from $C_{prd}^{\delta/2, \delta}({\mathbf Q})  $ to
$ C^{1+\delta/2, 2+\delta}({\mathbf Q})   $.
By periodicity of $m(0,\cdot)$
and next uniqueness of the solution $q$ in Lemma \ref{lemma:aux}, we further have $ m\in C_{prd}^{1+\delta/2, 2+\delta}({\mathbf Q})  $.

Next,  taking any ${\mathbb R}^n$-valued function
$ m\in  C_{prd}^{1+\delta/2, 2+\delta}({\mathbf Q})   $  in \eqref{pde4S}, we find a unique solution
$ S$ by solving a linear ODE for each  $\alpha\in {\mathbb R}^\iota$  and  define the mapping
\begin{align}\label{deflam2}
 S=\Lambda_2 ( m).
\end{align}

Let $\Phi(t,\tau)$ be the fundamental solution matrix of the ODE $\dot z(t)= -(A^T-P(t)BR^{-1} B^T)z(t)$.
Denote $|f|_0=|f|_{0;[0,T]}= \sup_{t\in [0,T]}|f(t)|.$  Define the constant
\begin{align}
K_0=&(T+T^{1-\delta/2}+ T^{2-\delta/2})\sup_{t,\tau\le T}\big[|\Phi(t,\tau)|+ |\partial_t\Phi(t, \tau)|\big]\notag   \\
&\times (|Q\Gamma_1-PD|_0+\kappa |P|_0 +| Q_f \Gamma_{1f} |).\notag
\end{align}
The estimate in the specific form in the next lemma will be useful for identifying a contraction condition later for a fixed point problem.

\begin{lemma}\label{lem:SK0}
For ${\mathbf Q}=[0,T]\times {\mathbb R}^\iota$ and $m(0,\cdot)\in C_{prd}^{2+\delta}({\mathbb R}^\iota)$,
$\Lambda_2$ in \eqref{deflam2} is a mapping from $ C_{prd}^{1+\delta/2, 2+\delta}({\mathbf Q})  $
to $  C_{prd}^{\delta/2, \delta}({\mathbf Q})  $ with
\begin{align}
|S|_{\delta/2, \delta; {\mathbf Q}}&\le K_0|m|_{1+\delta/2,2+\delta;{\mathbf Q}} \notag\\
&+  \sup_t |\Phi (t, T) Q_f \Gamma_{1f}|\cdot\Big( |m(0,\cdot)|_{0;{\mathbb R}^\iota}+ \sum_i |\partial_{\alpha_i} m(0,\cdot)|_{0;{\mathbb R}^\iota} \Big) \notag \\
&  + T^{1-\delta/2} \sup_t |\partial_t \Phi(t,T)|
 \cdot |Q_f\Gamma_{2f}\eta|
  +\sup_{t,\tau} |\Phi(t,\tau)|\cdot \big(T |Q\Gamma_2 \eta| + |Q_f \Gamma_{2f}\eta|\big)   .
 \label{sHolder}
\end{align}
\end{lemma}

\begin{proof}
See Appendix A.
\end{proof}

Now we introduce the following fixed point equation
\begin{align}
S= \Lambda_2 \Lambda_1 (S), \label{FPs}
\end{align}
which determines a solution of the system
\eqref{pde4S}--\eqref{pde4m}.

Consider $m=\Lambda_1  (S)$ and $ m'= \Lambda_1 (S')$,  with the same initial condition  $m(0,\alpha)=m'(0,\alpha)$, where  $S$ and $S'$ are any functions in
$  C_{prd}^{\delta/2, \delta}({\mathbf Q})  $.
By Lemma \ref{lemma:aux}, for some fixed constant $C_0$,
we have
\begin{align}
|m-m'|_{1+\delta/2, 2+\delta; {\mathbf Q}}\le C_0 |BR^{-1}B^T|\cdot | S-S'|_{\delta/2, \delta; {\mathbf Q}}.\label{mincrem}
\end{align}

\begin{theorem}\label{theorem:Sm}
Given $m(0,\cdot)\in C_{prd}^{2+\delta}({\mathbb R}^\iota)$,
the  PDE system \eqref{pde4S}--\eqref{pde4m}  for
$(S,m)$  has a unique solution in
$ C_{prd}^{\delta/2, \delta}({\mathbf Q})\times C_{prd}^{1+\delta/2, 2+\delta}({\mathbf Q})  $  if $C_0 K_0 |BR^{-1} B^T |<1,$
where $C_0$ is given in Lemma \ref{lemma:aux} taking ${\mathbf Q}=[0,T]\times {\mathbb R}^\iota$ and $F(t)=A-BR^{-1}B^T P(t)+D$.
\end{theorem}

\proof
For any two functions
$m, m'$ in  $C_{prd}^{1+\delta/2, 2+\delta}({\mathbf Q})$ with the same initial condition $m(0,\cdot)$,
 define $S=\Lambda_2 (m)$, $S'=\Lambda_2(m')$. Then we
follow the method of establishing \eqref{sHolder} to similarly obtain
\begin{align}
|S-S'|_{\delta/2, \delta; {\mathbf Q}}&=|\Lambda_2(m)-\Lambda_2(m')|_{\delta/2, \delta ; {\mathbf Q}}\le K_0 |m-m'|_{1+\delta/2, 2+\delta; {\mathbf Q}  }, %\label{ssmmd}
\end{align}
which together with \eqref{mincrem} yields
\begin{align}
|\Lambda_2 \Lambda_1 (S)- \Lambda_2\Lambda_1 (S') |_{\delta/2, \delta; {\mathbf Q} }
&\le K_0| \Lambda_1(S) -\Lambda_1(S') |_{1+\delta/2, 2+\delta ; {\mathbf Q}  }\notag \\
& \le C_0K_0|BR^{-1} B^T|\cdot |S-S'|_{\delta/2, \delta ; {\mathbf Q} }.\notag
\end{align}
We then apply a fixed point theorem on the Banach space $ C_{prd}^{\delta/2, \delta}({\mathbf Q}) $. \qed

\begin{remark}
The contraction condition is satisfied when $T$ is sufficiently small or the system has weak coupling, i.e., $|\Gamma_1|+|\Gamma_{1f}|+|D|+\kappa$ is sufficiently small.
\end{remark}

\begin{corollary}\label{coro:C4}
If $m(0,\cdot)\in C_{prd}^{2+\delta}({\mathbb R}^\iota)$ in Theorem \ref{theorem:Sm} is replaced by $m(0, \cdot)\in C_{prd}^{4+\delta}(\mathbb{R}^\iota)$, then we have $S\in C^{0,2}({\mathbf Q})$ and $m\in C^{1,4}({\mathbf Q})$, and moreover, $S$ (resp., $m$) has all up to second order partial derivatives with respect to $\alpha$, all belonging to $  C_{prd}^{\delta/2, \delta}({\mathbf Q}) $ (resp., $ C_{prd}^{1+\delta/2, 2+\delta}({\mathbf Q})$).
\end{corollary}

\proof Consider the case of $\iota=1$ while the case of $\iota=2$ is similar. We apply a difference quotient approach \cite[pp 72--74]{F08} by setting
$S_h(t,\alpha)=h^{-1}[S(t, \alpha+h)-S(t,\alpha)]$ and $m_h(t,\alpha)=h^{-1}[m(t, \alpha+h)-m(t,\alpha)] $, $h>0$, to obtain
 \begin{align}
\partial_t{S}_h(t,\alpha)&= -(A^T  -P(t)BR^{-1}B^T)S_h(t,\alpha)  \label{pde4Sh}  \\
&  \quad -\kappa P(t)  \Delta m_h(t,\alpha) +(Q\Gamma_1-P(t)D)  m_h(t,\alpha), \notag \\
\partial_t{m_h (t,\alpha)} &=\kappa\Delta m_h(t,\alpha) + (A-BR^{-1}B^T P(t)+D )m_h(t,\alpha )\label{pde4mh}  \\
&\quad  - B R^{-1} B^T S_h(t,\alpha), \qquad  0<t<T, \quad \kappa>0,    \notag
\end{align}
where  $S_h(T,\alpha)= -Q_f \Gamma_{1f}m_h(T,\alpha)$, with initial  condition $m_h(0,\cdot)$.
It is easy to show  $\sup_h |m_h(0,\cdot)|_{2+\delta;{\mathbb R}^\iota}\le C$.
Under the contraction condition of Theorem \ref{theorem:Sm}, $(S_h, m_h)$ is the unique solution of \eqref{pde4Sh}--\eqref{pde4mh}. In view of Lemma \ref{lemma:aux} and \eqref{sHolder}, we can show
$$\sup_h \big(|S_h|_{\delta/2, \delta;{\mathbf Q}}  +|m_h|_{1+\delta/2, 2+\delta; {\mathbf Q}}\big)\le C.
$$
By a standard subsequence argument \cite[pp 72--74]{F08}, we see that $(\partial_\alpha S   , \partial_\alpha m)$ is a solution of \eqref{pde4Sh}--\eqref{pde4mh} with initial condition $\partial_\alpha m(0, \alpha)$.
Repeating this argument, we further obtain $(\partial^2_\alpha S   , \partial^2_\alpha m)$ as a solution to   \eqref{pde4Sh}--\eqref{pde4mh}  with initial condition $\partial_\alpha^2 m(0, \alpha)$, so that $(\partial^2_\alpha S   , \partial^2_\alpha m)\in    C_{prd}^{\delta/2, \delta}({\mathbf Q}) \times C_{prd}^{1+\delta/2, 2+\delta}({\mathbf Q})$. \qed

\begin{figure}[t]
\begin{center}
\begin{tabular}{cc}
\psfig{file=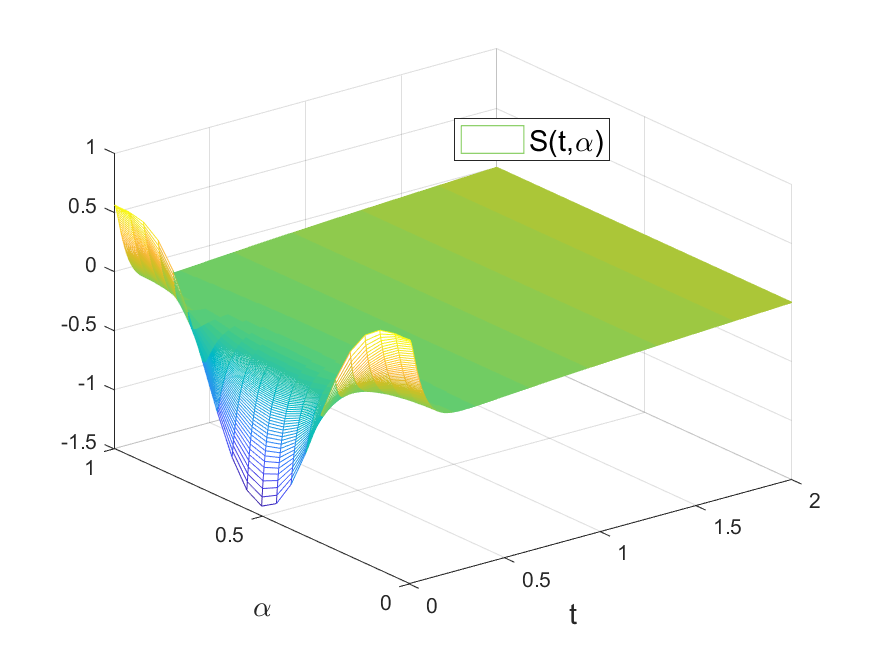, width=6cm, height=4cm}&
\psfig{file=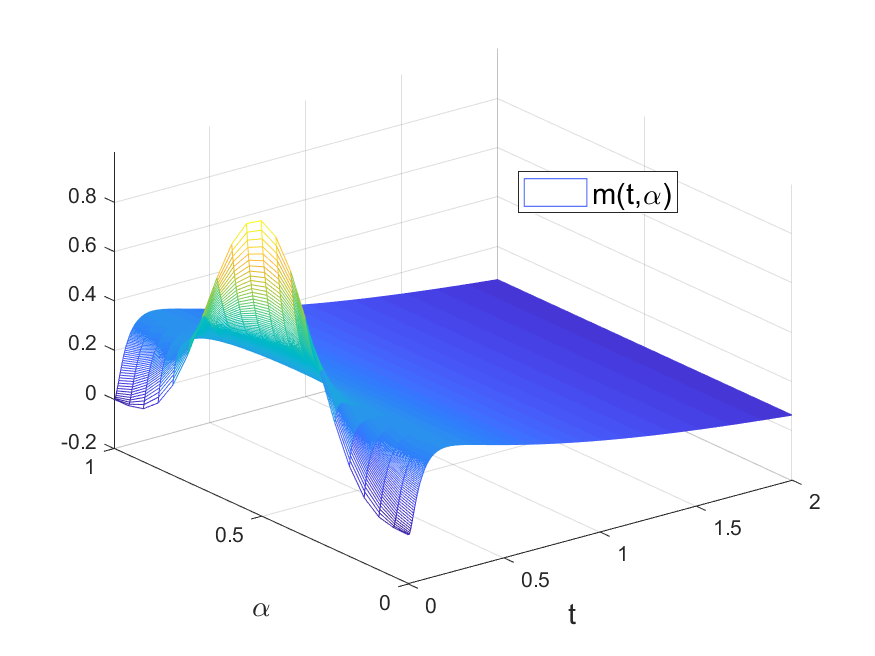, width=6cm, height=4cm}
\end{tabular}
\end{center}
\caption{ $S(t,\alpha)$ and $m(t,\alpha)$ in  \eqref{pde4S}--\eqref{pde4m} computed by iterations with
$t\in [0, 2]$ and
$\alpha \in [0, 1]$. }
\label{fig:fbpde}
\end{figure}

\subsection{Numerical example}
\begin{example}
 In the LQ model \eqref{dxat3}--\eqref{jakiu}  we consider a ring network with   $A=0.5$,  $B=1$,  $D=0.2$, $\kappa=0.5$, $Q=Q_f=2$,  $R=1$,   $\Gamma_1=0.8$,
$\Gamma_2=0.2$, $ \eta =1$, $\Gamma_{1f}=0.6$, $ \Gamma_{2f}=0.4$,
$T=2$,   and $m(0,\alpha)=\psi(\alpha):=256\alpha^4(1-\alpha)^4$,
 $0\le \alpha\le 1$.
\end{example}

If we extend  $\psi$ on $\mathbb{R}$ with period 1, then it can be shown that $\psi\in C_{prd}^{4+1/2}({\mathbb R})$.
The  system \eqref{pde4S}--\eqref{pde4m} is solved by a difference scheme, with stepwise $\Delta t=0.0025$ and
$\Delta \alpha =0.05$.  The numerical solution of $(S,m)$ is displayed in Fig. \ref{fig:fbpde}.

\section{Closed-loop analysis of the LQ model}

\label{sec:ne}

In this section,
 the cases of finite networks ${\mathbb T}_{1/N}$ and
 ${\mathbb T}^2_{1/N}  $ are considered.
 Let the LQ model be specified by \eqref{dxat3}--\eqref{jakiu}.
For tractability of the subsequent estimates,  all remaining sections  consider  scalar states, control and Brownian motions in \eqref{dxat3}--\eqref{jakiu} so that
$n=n_1=n_2=1$.

We introduce the following  assumptions (H1)--(H4), which are in effect throughout Sections \ref{sec:ne}, \ref{sec:so} and \ref{sec:eNE}.

(H1) For a given function $m_0(\alpha)$ from $[0,1)^\iota$ to
${\mathbb R}$, we have
$x^i_{\alpha_k^N}(0)=m_0({\alpha_k^N})$, {for all} $ i,k$.
The periodic extension of $m_0(\cdot)$ to the domain ${\mathbb R}^\iota$, denoted by $\widetilde m_0$, belongs to
$ C^{4+\delta}_{prd}( {\mathbb R}^\iota)$ for some
$\delta\in (0,1)$.

(H2) There exists a unique solution $( S,m) $ to  \eqref{pde4S}--\eqref{pde4m}
with $n=1$ on ${\mathbf Q}=[0,T]\times {\mathbb R}^\iota$, such that
$S\in C_{prd}^{\delta/2, \delta}(\mathbf{Q})$, $m\in C_{prd}^{1+\delta/2, 2+\delta}(\mathbf{Q})$,
 $m(0,\alpha)=\widetilde m_0(\alpha) $.

(H3) The solution $(S,m) $ to \eqref{pde4S}--\eqref{pde4m} is such that  the 2nd (resp., the 4th) order partial derivatives of $  S(t,\alpha)$ (resp., of $ m(t,\alpha)$)  with respect to $\alpha$
 are  continuous in $(t,\alpha)\in  [0,T]\times {\mathbb R}^\iota$.

(H4) $\lim_{N\to \infty} N_1/N^3=\infty$  for the  ring network  ${\mathbb T}_{1/N}$;  $\lim_{N\to\infty} N_1/(N^2\ln N)=\infty$  for the  torus network ${\mathbb T}_{1/N}^2$.
\begin{remark}
(i) We assume  equal initial states within the same cluster to streamline the presentation, and can deal with random initial states easily.
(ii) For  higher order regularity of the solution in (H3), Corollary \ref{coro:C4} provides sufficient conditions.
(iii) Since the second order interaction terms in
\eqref{dxat3} are random and scaled up by $N^2$, we need to average with fast growth of $N_1$ as in (H4) to reduce their  fluctuations.
\end{remark}

By assigning the agents' actual nodal positions to $\alpha$ in  \eqref{bru}, the LMFG-based decentralized strategies for the finite population take the form
\begin{align}
\hat u^i_{\alpha_k^N}(t)=& -R^{-1} B[P(t)x^i_{\alpha_k^N}(t)+S(t,{\alpha_k^N}   )],  \quad  1\le i\le N_1, \quad 1\le k\le N^*. \label{uhati}
\end{align}
In further analysis, for notational brevity, we will often drop the time argument in a process. For instance, we simply write $x^i_{\alpha_k^N}(t) $ as $x^i_{\alpha_k^N}$.
The strategies in \eqref{uhati}  yield the closed-loop dynamics:
\begin{align}
{\rm System\ (A):}&\notag \\
d \hat x^i_{\alpha_k^N}=\ &\{ [A  - B^2 R^{-1}  P(t)] \hat x_{\alpha_k^N}^i   +D
 \hat x^{(N_1)}_{\alpha_k^N}    - B^2R^{-1} S(t,{\alpha_k^N})\}  dt   \notag
\\ &\hskip -0.5cm +\kappa N^2 \sum_{k':\alpha_{k'}^N\sim \alpha_{k}^N   }\big(  \hat x^{(N_1)}_{\alpha_{k'}^N}-  \hat x^{(N_1)}_{\alpha_k^N} \big)dt
   +\sigma dw^i_{\alpha_k^N} ,   \label{xiakn}\\
   &  1\le i\le N_1,  \quad 1\le k\le N^*, \notag
\end{align}
with  $ \hat x^i_{\alpha_k^N}(0)=    x^i_{\alpha_k^N} (0)  $ and
$ \hat x_{\alpha_k^N}^{(N_1)} = \frac{1}{N_1} \sum_{j=1}^{N_1} \hat x^j_{\alpha_k^N}.$

In the $\epsilon$-Nash equilibrium analysis, only one agent applies a different strategy while all the other agents stick to the strategies in \eqref{uhati}.
Without loss of generality, we consider the strategy change of  agent $ {\cal A}^1_{\alpha_1^N}$ within the cluster residing  at node $\alpha_1^N$,  and the system dynamics are given by
\begin{align}
{\rm System\ (B):}&\notag \\
dx^1_{\alpha_1^N} =\ &[ Ax^1_{\alpha_1^N}
+ Bu^1_{\alpha_1^N}
+ D  x_{\alpha_1^N}^{(N_1)}]dt\notag \\
&+
\kappa N^2 \sum_{k:\alpha_k^N\sim \alpha_1^N } (x_{\alpha_k^N}^{(N_1)}
-x_{\alpha_1^N}^{(N_1)})dt +\sigma dw_{\alpha_1^N}^1,\notag\\
dx^i_{\alpha_k^N} =\ &[ Ax^i_{\alpha_k^N}
+ B\hat u^i_{\alpha_k^N}
+ D  x_{\alpha_k^N}^{(N_1)}]dt\notag \\
&+
\kappa N^2  \sum_{k': \alpha_{k'}^N \sim \alpha_{k}^N  } (x_{\alpha_{k'}^N}^{(N_1)}
-x_{\alpha_k^N}^{(N_1)})dt+\sigma dw_{\alpha_k^N}^i,\notag
\end{align}
for $i\in\{1, \cdots, N_1\}$ if $k\ne 1$; and $i\in \{2, \cdots, N_1\}$ if $k=1$.

The next theorem gives the a priori upper bound of the cost that is achieved by the decentralized strategies in \eqref{uhati}.
\begin{theorem}
\label{prop:Jbound}
There exists a constant $\hat C_0$  such that
$$
\sup_{N,N_1}\max_{k\le N^*, i\le N_1}J_{\alpha_k^N}^i(\hat u^i_{\alpha_k^N},
\hat u^{-i}_{\alpha_k^N},
\hat {\bf u}_{\alpha_{-k}^N})\le \hat C_0.
$$
\end{theorem}

We postpone its proof, which is technical and needs a series of preparations.
For the $\epsilon$-Nash equilibrium analysis, Theorem \ref{prop:Jbound} provides a benchmark performance level when one agent  unilaterally  applies  an alternative strategy  $u^i_{\alpha_k^N}  $ for its own improvement.
 In this case it suffices to consider  $ u^i_{\alpha_k^N}$ satisfying
 \begin{align}
  J_{\alpha_k^N}^i(u^i_{\alpha_k^N},\hat u^{-i}_{\alpha_k^N}, \hat {\bf u}_{\alpha_{-k}^N}  ) \le \hat C_0. \label{pcbc0}
\end{align}
The restriction to such  strategies $u^i_{\alpha_k^N}$ will allow us to obtain uniform error bounds when approximating  $x^{(N_1)}_{\alpha_k^N}$ and
$N^2  \sum_{k': \alpha_{k'}^N \sim \alpha_{k}^N  } (x_{\alpha_{k'}^N}^{(N_1)}
-x_{\alpha_k^N}^{(N_1)}) $, $1\le k\le N^*$, generated by $( u^i_{\alpha_k^N},\hat u^{-i}_{\alpha_k^N}, \hat {\bf u}_{\alpha_{-k}^N}    )$.

 By \eqref{xiakn}, we have
\begin{align}
d \hat x^{(N_1)}_{\alpha_k^N}=& [A-B^2R^{-1} P(t)+D] \hat x^{(N_1)}_{\alpha_k^N} dt
- B^2R^{-1} S(t,{\alpha_k^N}) dt \notag    \\
&+ \kappa N^2 \sum_{k':\alpha_{k'}^N\sim  \alpha_{k}^N   }(\hat x^{(N_1)}_{\alpha_{k'}^N}
- \hat x^{(N_1)}_{\alpha_k^N} ) dt  +\sigma d w^{(N_1)}_{\alpha_k^N}  , \label{xkn4}
\end{align}
where
$w^{(N_1)}_{\alpha_k^N}=\frac{1}{N_1} \sum_{j=1}^{N_1}
w^j_{\alpha_k^N}.$

To find the cost bound in Theorem \ref{prop:Jbound}, a  key step is to estimate the following quantities in System (A):
\begin{align}\label{exxn2}
\sup_{N,N_1}\sup_{k,t}\EE|\hat x^{(N_1)}_{\alpha_k^N} (t)  |^2, \qquad
\sup_{N,N_1}\sup_{k,t}N^2\EE\Big| \sum_{k':\alpha_{k'}^N\sim  \alpha_{k}^N   }(\hat x^{(N_1)}_{\alpha_{k'}^N}
- \hat x^{(N_1)}_{\alpha_k^N} )(t)\Big|^2.
\end{align}
As it turns out below, the resulting analysis is much more difficult than in GMFGs with dense networks, which is due to the second order interaction term in \eqref{dxat3}.

To facilitate the estimate of  $\hat x^{(N_1)}_{\alpha_k^N}$, we introduce an auxiliary  deterministic system:
\begin{align}
{\rm System\ (C)}:& \notag \\
\frac{d}{dt} \bar x_{\alpha_k^N} = & [A-B^2R^{-1} P(t)+D] \bar x_{\alpha_k^N}
- B^2R^{-1} S(t,{\alpha_k^N})  \notag    \\
&+ \kappa N^2 \sum_{k': \alpha_{k'}^N \sim \alpha_{k}^N  }(\bar x_{\alpha_{k'}^N}
-\bar x_{\alpha_k^N} ) ,  \quad  1\le k\le N^*,  \label{checkxd}
\end{align}
where  $\bar x_{\alpha_k^N}(0)= m_0({\alpha_k^N})$.

\begin{lemma}\label{lem:xbBound}
There exists a constant $C$ such that
$\sup_{N,k\le N^*,t\le T}|\bar x_{\alpha_k^N}(t)|\le C$.
\end{lemma}

\begin{proof}
For $t\in [0,T]$, define $x^N_{\rm max}(t)=\max_{1\le k\le N^*}\bar x_{\alpha_k^N}(t)$ which is then absolutely continuous and so differentiable almost everywhere on $[0,T]$. If $ x_{\rm max}^N(t)$ is differentiable at $t$, denote the set ${\mathbf s}_t^N=\{k: x_{\rm max}^N(t)= \bar x_{\alpha_k^N}(t)\}$.  Consider the sequence $\{x_{\rm max}^N(t+1/j), j\ge j_0 \}$ for a sufficiently large $j_0$ if $t<T$. There necessarily exists some $\hat k\in {\mathbf s}_t^N$  such that  $x_{\rm max}^N(t+1/j)= \bar x_{\alpha_{\hat k}^N}(t+1/j)$ for an infinite number of values of $j$  (we replace $t+1/j$ by $t-1/j$ if $t=T$). Therefore, we have $\frac{d}{dt}x^N_{\rm max}(t) =
\frac{d}{dt}\bar x_{\alpha_{\hat k}^N}(t)   $, which further implies that
whenever $x^N_{\rm max}(t)$ is positive, its derivative, if existing, is at most  $|A-B^2R^{-1} P(t)+D| x_{\rm max}^N (t)+c_1$ for
$c_1\coloneqq B^2R^{-1} \sup_{\alpha,t} |  S(t,\alpha)|$ since $\kappa >0$.
 To bound $ x_{\rm max}^N(t)$ on $[0,T]$ from above,
 we construct the auxiliary ODE:
 \begin{align}\label{yc1}
 \dot y(t) = |A-B^2R^{-1} P(t)+D| y (t)+c_1  , \quad y(0)= \max_{\alpha\in [0,1)^\iota} |m_0(\alpha)| .
 \end{align}
 Then by Lemma \ref{lem:odecom} (i), we have
 $ x_{\rm max}^N(t)\le y(t)$   for all $ t\in [0,T]$.

Similarly, we apply Lemma \ref{lem:odecom} (ii) to obtain a constant lower bound for $x^N_{\rm min}(t)\coloneqq\min_{1\le k\le N^*}\bar x_{\alpha_k^N}(t)$.  This completes the proof.
\end{proof}

\begin{lemma} \label{lem:xbm5}
We have
$\sup_{k\le N^*,t\le T}|\bar x_{\alpha_k^N}(t)-m(t,\alpha_k^N) |=o(1)$ as $N\to \infty$.
\end{lemma}

\begin{proof}
Taking $\alpha=\alpha_k^N$ in \eqref{pde4m}, we write the equation as a
perturbed ODE:
\begin{align}
\frac{d}{dt} m(t,\alpha_k^N)&= [A-B^2R^{-1} P(t)+D] m(t,\alpha_k^N) -B^2 R^{-1} S(t,\alpha_k^N) \label{dtme}  \\
 & +  \kappa N^2 \sum_{k': \alpha_{k'}^N \sim \alpha_{k}^N  }[m (t,\alpha_{k'}^N)
-m(t,\alpha_k^N) ]  +\epsilon_{\alpha_k^N}(t), \notag
\end{align}
where
$$
\epsilon_{\alpha_k^N}(t)= \kappa \Delta m(t,\alpha_k^N) - \kappa N^2 \sum_{k': \alpha_{k'}^N \sim \alpha_{k}^N  }[m (t,\alpha_{k'}^N)
-m(t,\alpha_k^N) ]  ,
$$
and
$ \sup_{k,t} |\epsilon_{\alpha_k^N}(t)|\to 0 $ as $N\to \infty$ since $m\in  C^{1+\delta/2,2+\delta}_{prd}([0,T]\times {\mathbb R}^\iota) $.
Denote $\zeta_{\alpha_k^N}(t)= \bar x_{\alpha_k^N}(t)- m(t,\alpha_k^N)$. By  \eqref{checkxd} and  \eqref{dtme},  we have
\begin{align}
\frac{d}{dt} \zeta_{\alpha_k^N} (t)&= [A-B^2R^{-1} P(t)+D] \zeta_{\alpha_k^N}(t)
\notag \\
 & +  \kappa N^2 \sum_{k': \alpha_{k'}^N \sim \alpha_{k}^N  }[\zeta_{\alpha_{k'}^N}(t)
-\zeta_{\alpha_k^N}(t) ] -\epsilon_{\alpha_k^N}(t),\quad k\le N^*.\notag
\end{align}
Following the method in proving Lemma
\ref{lem:xbBound}, we construct the auxiliary ODE \begin{align}
\dot{y}(t)= |A-B^2R^{-1} P(t)| y(t)+ \sup_{t,k} |\epsilon_{\alpha_k^N}(t)|, \quad y(0)=0, \notag
\end{align}
and   obtain an upper bound of the form
$
\max_{t,k}\zeta_{\alpha_k^N} (t)\le C \sup_{t,k} |\epsilon_{\alpha_k^N}(t)|.
$
Similarly we obtain a lower bound for $ \min_{t,k}\zeta_{\alpha_k^N} (t) $, and conclude that
$ \max_{t,k}|\zeta_{\alpha_k^N} (t)|\le C \sup_{t,k} |\epsilon_{\alpha_k^N}(t)|
=o(1)$ as $N\to \infty$.
\end{proof}

Denote
\begin{align}\label{checkx}
{\check x}_{\alpha_k^N}(t)= \hat x^{(N_1)}_{\alpha_k^N}(t)
- \bar x_{\alpha_k^N}(t), \qquad  k\le N^*,\quad  0\le t\le T.
\end{align}

To obtain bounds for the two terms in \eqref{exxn2},
we start by analyzing ${\check x}_{\alpha_k^N}(t)$, $1\le k\le N^*$, which will be stacked into a large vector by the following rule:
\begin{enumerate}
\item[(i)] For the ring network ${\mathbb T}_{1/N}$,   define
\begin{align}\label{xwN}
{{\check {\bf x}}}= [{\check x}_{\alpha_1^N}, \cdots ,{\check x}_{\alpha_N^N}]^T , \quad
{\bf  w}= [ w^{(N_1)}_{\alpha_1^N},\cdots,  w^{(N_1)}_{\alpha_N^N} ]^T,
\end{align}
where  $\alpha_k^N= \frac{k-1}{N} +\frac{1}{2N}$, $1\le k\le N$ (see Section  \ref{sec:rtop}).

\item[(ii)] For the torus network ${\mathbb T}^2_{1/N}$,  let the
$N^2$ nodes  be arranged into an
$N\times N$ array. The entry at the $i$th row and the $j$th column, to be called the $(i,j)$-th node, is
$(\frac{i-1}{N} +\frac{1}{2N}, \frac{j-1}{N} +\frac{1}{2N} )$, $1\le i,j\le N$. Each node has 4 neighbors.  For instance, the $(1,1)$-th node has neighbors: the  $(1,2)$-th node, the $(1,N)$-th node, the $(2,1)$-th node, and the $(N,1)$-th node.
We list the $N^2$ nodes into a single column, starting with the first row,  next the second row, and so forth until the last row.
Following the above order,
we now list the $N^2$ terms $\check x_{\alpha_k^N}$ (resp., $w^{(N_1)}_{\alpha_k^N} $) into a  large vector denoted by $\check{\bf  x}$ (resp., ${\bf w}$):
\begin{align}\label{xwN2}
\check {\bf x}= [\check x_{\alpha_1^N}, \cdots, \check x_{\alpha_{N^2}^N}]^T,
\quad  {\bf w}= [w_{\alpha_1^N}^{(N_1)}, \cdots, w^{(N_1)}_{\alpha_{N^2}^N}]^T,
\end{align}
where $\alpha_1^N$ and $\alpha_{N^2}^N$ are the $(1,1)$-th and the $(N,N)$-th nodes, respectively.
\end{enumerate}

Our subsequent analysis  will give details for the torus case, and the ring case can be similarly handled and is simpler.
For the error estimate, we will derive an equation for $\check{\bf x}$.
As a preparation,
denote the $N\times N$ circulant matrix \cite{D79}
\begin{align*}
{\cal R}_N=
\begin{bmatrix}
0 & 1 & 0 & \cdots &0&1 \\
1 & 0 & 1 & \cdots &0&0\\
 0 & 1 & 0 & \cdots &0&0 \\
 0 & 0 & 1 & \cdots & 0 &0\\
 \vdots &\vdots &\vdots &\ddots &\vdots & \vdots \\
 1 & 0 & 0 & \cdots &1&0
\end{bmatrix},
\end{align*}
which is symmetric and is the Laplacian matrix of the ring network ${\mathbb T}_{1/N}$.

 For the ring network ${\mathbb T}_{1/N}$, letting
   $\hat A(t)=A-B^2R^{-1}P(t)+D$ and
 $\hat A_N(t)=\hat A(t)-2\kappa N^2$,
we introduce the following circulant matrix
\begin{align}
{\mathcal M}^{rn}_N(t)&=
\begin{bmatrix}
\hat A_N & \kappa N^2 & 0& \cdots& 0& \kappa N^2   \\
\kappa N^2 & \hat A_N & \kappa N^2 & \cdots & 0 & 0\\
\vdots & \vdots & \vdots & \vdots & \vdots  & \vdots \\
 \kappa  N^2&0 &0 & \cdots & \kappa N^2 &  \hat A_N
\end{bmatrix} \label{mnrnt}   \\
&= \hat A(t) I_N-2\kappa N^2 I_N +\kappa N^2{\cal R}_N. \notag
\end{align}
where each row has 3 non-zero entries.

To analyze $\check {\bf x}$ for the torus network ${\mathbb T}_{1/N}^2$,  define
\begin{align}\label{rksum}
{\cal T}_N=  I_{N}\otimes {\cal R}_N +{\cal R}_N\otimes I_{N},
\end{align}
which is the Laplacian of  ${\mathbb T}^2_{1/N}$.
An example of ${\cal T}_N$ with $N=4$ is illustrated in Appendix B.
Following the definition in \cite{HJ91}, the right hand side of \eqref{rksum} is the Kronecker sum $ {\cal R}_N\oplus {\cal R}_N$. We further denote
\begin{align}\label{mntrn}
{\cal M}^{tr}_N(t)=I_{N^2}\otimes(\hat A(t) -4\kappa N^2) + \kappa N^2{\cal T}_N.
\end{align}
 So ${\cal M}_N^{tr}$ is a symmetric matrix with  dimension
$N^2 \times N^2 $.

We have the following representation of the process $\check{\bf x}$:
\begin{lemma}\label{lemma:xmw}
Let $\check{x}_{\alpha_k^N}(t)$ be defined by \eqref{checkx}. Then
\begin{itemize}
\item[(i)] for the ring network ${\mathbb T}_{1/N}$ and
\eqref{xwN}, we have
\begin{align}
d\check {\bf x}(t) = {\cal M}_N^{rn} \check{\bf x}(t) dt +\sigma d{\bf w}(t),
\end{align}
with initial condition $\check{x}_{\alpha_k^N}(0)$, $1\le k\le N$;
\item[(ii)] for the torus network  ${\mathbb T}^2_{1/N}  $ and \eqref{xwN2}, we have
 \begin{align}\label{xcms5}
d\check {\bf x}(t) = {\cal M}_N^{tr} \check{\bf x}(t) dt +\sigma d{\bf w}(t),
\end{align}
with initial condition $\check{x}_{\alpha_k^N}(0)$, $1\le k\le N^2$.
 \end{itemize}
\end{lemma}

\begin{proof}
By \eqref{xkn4} and \eqref{checkxd}, it is straightforward to show
\begin{align}%\label{sdexch}
d{\check x}_{\alpha_k^N}&  =  [A-B^2R^{-1} P(t)+D]{\check x}_{\alpha_k^N} dt+ \kappa N^2\sum_{k': \alpha_{k'}^N\sim  \alpha_{k}^N  } ({\check x}_{\alpha_{k'}^N}  -{\check x}_{\alpha_k^N}   ) dt \label{xw42}  \\
& + \sigma d w^{(N_1)}_{\alpha_k^N}, \quad 1\le k\le N^*, \notag
\end{align}
where ${\check x}_{\alpha_k^N}(0) =0$.
Arranging ${\check x}_{\alpha_k^N}$ into a vector as in \eqref{xwN}--\eqref{xwN2}, and  checking the coefficient matrix of each case, we derive the  SDEs in (i) and (ii).
\end{proof}

Define the polynomial $\varphi(z)=z+z^{N-1}$ for ${\cal R}_N$. By the results on general circulant  matrices (see \cite{D79}), the eigenvalues of ${\cal R}_N$ are given by
\begin{align}\label{lamkN}
\lambda_k=\varphi (e^{{\bf i} 2\pi (k-1)/N}) = 2\cos \frac{2\pi(k-1)}{N}  , \quad  {\bf i}^2=-1,\  1\le k\le N.
\end{align}
Denote $\omega_k= e^{{\bf i} 2\pi (k-1)/N} $, $1\le k\le N$, and the Fourier matrix
$$
F=\frac{1}{\sqrt{N}}
\begin{bmatrix}
1 & 1 & 1 & \cdots & 1 \\
1 & \omega_N  &\omega_N^2    &\cdots       & \omega_{N}^{N-1} \\
\vdots & \vdots &\vdots &\vdots & \vdots \\
1 & \omega_N^{N-1}  &\omega_N^{2(N-1)}   &\cdots   & \omega_N^{(N-1)^2}
\end{bmatrix}.
$$
By  results on circulant matrices \cite{D79},  $F$ is unitary (i.e., $F^*F=I$) and
\begin{align}\label{RNdg}
F^*{\cal  R}_N F= \mbox{diag}[\lambda_1, \cdots, \lambda_N  ],
\end{align}
where $F^*$ is the  conjugate transpose of $F$.
Denote
$$
\lambda_{jk} =\lambda_j+\lambda_k = 2 \left(\cos \frac{2\pi(j-1)}{N} + \cos \frac{2\pi(k-1)}{N} \right) ,\quad 1\le j,k\le N,
$$
and the $N^2\times N^2$ matrix
\begin{align*}
\Lambda_0&=\mbox{diag}[\lambda_{11}, \cdots, \lambda_{1N}, \lambda_{21}, \cdots, \lambda_{2N},\cdots, \lambda_{N1}, \cdots, \lambda_{NN}].
\end{align*}

\begin{lemma} \label{lem:FMd}
The following properties hold:
\begin{itemize}

\item[(i)]  ${\cal M}^{rn}_N(t)$ in \eqref{mnrnt} is symmetric and
${\cal M}^{rn}_N(t){\cal M}^{rn}_N(s)={\cal  M}^{rn}_N(s) {\cal M}^{rn}_N(t)$ for all $t,s\in [0,T]$;

\item[(ii)]  ${\cal T}_N$ in \eqref{rksum} is symmetric;

\item[(iii)]  ${\cal M}^{tr}_N(t)$ in \eqref{mntrn} is symmetric and
 ${\cal M}^{tr}_N(t){\cal M}^{tr}_N(s)={\cal  M}^{tr}_N(s) {\cal M}^{tr}_N(t)$ for all $t,s\in [0,T]$;

\item[(iv)]
$(F\otimes F)^*{\cal M}^{tr}_N(t) (F\otimes F)= (\hat A(t) -4\kappa N^2) I_{N^2}+
\kappa N^2 \Lambda_0$,
and $F\otimes F$ is a unitary matrix.
\end{itemize}
\end{lemma}

\begin{proof}
(i) holds since ${\cal R}_N$ is symmetric. Using the property $(M_1\otimes M_2)^T= M_1^T\otimes M_2^T$ for any two matrices $M_1$ and $M_2$, we obtain (ii), which further implies (iii).
By basic properties of Kronecker products, it is easily shown that $F\otimes F$ is a unitary matrix.
By \cite[Theorem 4.4.5]{HJ91},
the $N^2$ eigenvalues $\lambda_{jk}$ (including algebraic multiplicities) of
${\cal T}_N$ are given by $\Lambda_0$, and we use properties of Kronecker products (see \cite[Lemma 4.2.10]{HJ91}) and \eqref{RNdg} to calculate
 \begin{align*}
(F\otimes F)^* {\cal  T}_N (F\otimes F)&=(F^*\otimes F^*)(I_N\otimes {\cal R}_N) (F\otimes F)\\
&\quad  +(F^*\otimes F^*)( {\cal R}_N \otimes I_N) (F\otimes F) \\
&  = (F^*F) \otimes (F^* {\cal R}_NF)+(F^* {\cal R}_NF)\otimes (F^*F)   \\
&=   \Lambda_0,
\end{align*}
which leads to (iv).
\end{proof}

By Lemma \ref{lem:FMd} (iv),
 all $N^2$ eigenvalues of ${\cal M}_N^{tr}(t)$ are given by
\begin{align}
\hat A(t) -2\kappa N^2 \left(2-\cos \frac{2\pi(j-1)}{N} - \cos \frac{2\pi(k-1)}{N}\right),
\quad 1\le j,k\le N.
\end{align}

Consider the ODEs
\begin{align}
\dot{\bf x}= {\cal M}_N^{rn} (t) {\bf x},\qquad
\dot{\bf x}= {\cal M}_N^{tr} (t) {\bf x},
\end{align}
where the dimension of ${\bf x}$ is compatible with the coefficient matrix.
In view of  Lemma \ref{lem:FMd} (i) and  (iii), we  apply
 Lemma \ref{lem:odeF} to get  the fundamental solution matrices
\begin{align}\label{Ph12ts}
\Phi^{ rn}(t,s)= e^{\int_s^t {\cal M}_N^{rn}(\tau) d\tau},\quad
\Phi^{ tr}(t,s)= e^{\int_s^t {\cal M}_N^{tr}(\tau) d\tau}, \quad t,s\in [0,T],
\end{align}
which are symmetric for all given $t,s\in [0,T]$.
We give some estimate for
\eqref{checkx}.

 \begin{lemma} \label{lem:xchBound}
{\rm (i)} For the ring network ${\mathbb T}_{1/N}$, we have
$$
\sup_{t,k} {\mathbb E}|{\check x}_{\alpha_k^N}(t)|^2 \le C/(NN_1).
$$

{\rm (ii)} For the torus network  ${\mathbb T}^{2}_{1/N}$,  we have
$$
\sup_{t,k} {\mathbb E}|{\check x}_{\alpha_k^N}(t)|^2 \le C \ln N/(N^2N_1).
$$

\end{lemma}

\begin{proof}
(i) For $\check {\bf x}$ given by \eqref{xwN},
we have
$\check {\bf x} (t)=\sigma \int_0^t \Phi^{rn} (t,\tau) d{\bf  w}({\tau})$
by Lemma \ref{lemma:xmw},
so that
\begin{align}
\EE | \check {\bf x}(t)|^2
&=\frac{\sigma^2}{N_1} {\rm Tr}\int_0^t e^{2  \int_\tau^t{\cal  M}_N^{rn}(r)dr} d\tau,\label{5xdtau}
\end{align}
by symmetry of ${\Phi}^{rn}(t,\tau)$.
By digitalization using \eqref{RNdg}, we have
 \begin{align}
 {\rm Tr}\int_0^t e^{2  \int_\tau^t{\cal  M}_N^{rn}(r)dr} d\tau
 &\le C\Big( 1+ \sum_{1\le k\le N-1} \frac{1}{N^2(1-\cos \frac{2\pi k}{N})}  \Big).
\end{align}
We estimate the sum
$ s_N\coloneqq\sum_{k=1}^{N-1} N^{-2} \left(1-\cos \frac{2\pi k}{N}\right)^{-1}.
$
There exists a fixed small $\theta_0>0$ such that $1-\cos \theta \ge \theta^2/4$ for all $|\theta|\le \theta_0$. We write
\begin{align}\label{5s123}
s_N= \sum_{k:0< 2\pi k/N\le \theta_0 } +  \sum_{k: \theta_0< 2\pi k/N<2\pi- \theta_0 }
+
\sum_{k: 2\pi-\theta_0\le  2\pi k/N< 2\pi }\eqcolon s_{N,1}+ s_{N,2} +s_{N,3},
\end{align}
where the summands are obvious and not displayed.
We obtain
\begin{align}
0\le s_{N,1}\le  \sum_{k: 0<  2\pi k/N\le  \theta_0 } 4 /(2\pi k)^2=O(1),
\end{align}
and  similarly obtain $s_{N,3}=O(1)$ while clearly $s_{N,2}=O(1)$. Hence, $s_N=O(1)$ and   $\EE | \check {\bf x}(t)|^2\le C/N_1$. By symmetry, we get the bound of each component.

 (ii) In analogue to \eqref{5xdtau},
we have
\begin{align}
\EE | \check {\bf x}(t)|^2
&=\frac{\sigma^2}{N_1} {\rm Tr}\int_0^t e^{2  \int_\tau^t{\cal  M}_N^{tr}(r)dr} d\tau. \notag
\end{align}
We have
\begin{align}
& {\rm Tr}\int_0^t e^{2  \int_\tau^t{\cal  M}_N^{tr}(r)dr} d\tau \label{trlog} \\
&={\rm Tr}\int_0^t e^{2 \int_\tau^t[ \kappa N^2 \Lambda_0 +(\hat A(r) -4\kappa N^2)I_{N^2} ]dr} d\tau
  \notag \\
& = \int_0^t \sum_{1\le i,j\le N} e^{2(t-\tau) (\lambda_i +\lambda_j -4)\kappa N^2}
e^{2\int_\tau^t \hat A(r)dr } d\tau \notag  \\
&\le CT+ C \sum_{i+j\ge 1, 0\le i,j\le N-1} \frac{1}{
N^2[ (1-\cos \frac{2\pi i}{N}) + (1-\cos \frac{2\pi j}{N})]
}  \notag  \\
&=O(\ln N), \notag
\end{align}
where we have obtained the last line using the method in \eqref{5s123} and a double integral estimate. Specifically, we estimate by splitting the last summation over $(i,j)$ into three parts for (i) $(i,j)$ with $i+j\ge 1$, $0\le i,j\le \theta_0 N/(2\pi)$, (ii) $(i,j)$ with $  N(1-\theta_0/(2\pi))\le i,j<N$, (iii) all pairs $(i,j)$ left out by (i) and (ii).
 Finally   for some constant $C$, we get
 ${\mathbb E}|{{\check {\bf x}}}(t)|^2\le C\ln N/{N_1}$.
By symmetry, we obtain the upper bound for each component.
\end{proof}

\begin{lemma}\label{prop:exhBound}
For some constant $C$, we have
$$
\sup_{N, N_1}\sup_{t,k}{\mathbb E}|\hat x^{(N_1)}_{\alpha_k^N}(t)|^2\le C.
$$
\end{lemma}
\proof  This lemma follows from \eqref{checkx}, Lemmas \ref{lem:xchBound}
and \ref{lem:xbBound}.
\qed

 \section{Estimate of the second order interaction term}
\label{sec:so}

This section analyzes the torus network ${\mathbb T}^2_{1/N}$ with detail.
For $\alpha_k^N$ in the array of $N^2$ nodes of ${\mathbb T}^2_{1/N}$ (see Sec. \ref{sec:ne} and \eqref{xwN2}), let $\alpha_{k_+}^N$ and $\alpha^N_{k_-}$ denote its two neighbors along the horizontal direction.
 Our next step is to approximate
 $
  \hat x^{(N_1)}_{\alpha_{k_+}^N}
+ \hat x^{(N_1)}_{\alpha_{k_-}^N}
-2 \hat x^{(N_1)}_{\alpha_k^N}
 $ contained in \eqref{xiakn} by the deterministic term $
 \bar x_{\alpha_{k_+}^N}
+\bar x_{\alpha_{k_-}^N}
-2\bar x_{\alpha_k^N} $.
We have the relation
\begin{align*}
&\hat x^{(N_1)}_{\alpha_{k_+}^N}
+ \hat x^{(N_1)}_{\alpha_{k_-}^N}
-2 \hat x^{(N_1)}_{\alpha_k^N} -( \bar x_{\alpha_{k_+}^N}
+\bar x_{\alpha_{k_-}^N}
-2\bar x_{\alpha_k^N})\\
& = {\check x}_{\alpha_{k_+}^N}+ {\check x}_{\alpha_{k_-}^N} -
2{\check x}_{\alpha_k^N},
\end{align*}
where ${\check x}_{\alpha_k^N} $ is defined by \eqref{checkx}.
Now we set
$v_{\alpha_k^N}= {\check x}_{\alpha_{k_+}^N}+ {\check x}_{\alpha_{k_-}^N} -
2{\check x}_{\alpha_k^N}$ and $w^{(N_1)}_{\alpha_{k}^N}=\frac{1}{N_1}\sum_{i=1}^{N_1}
w_{\alpha_k^N}^i $.
By  SDE \eqref{xw42} for $ \check x_{\alpha_k^N} $, we have
 \begin{align}
 dv_{\alpha_k^N}=& [A-B^2R^{-1} P(t)+D] v_{\alpha_k^N} dt
 + \kappa  N^2 \sum_{j: \alpha_j^N\sim \alpha_k^N}(v_{\alpha_j^N} -v_{\alpha_k^N}) dt \notag \\
  &+\sigma d w^{(N_1)}_{\alpha_{k_+}^N}+\sigma  d w^{(N_1)}_{\alpha_{k_-}^N}
  -2\sigma d w^{(N_1)}_{\alpha_{k}^N}, \qquad 1\le k\le N^2.\notag
 \end{align}
Denote
$ {\bf v}(t)= [v_{\alpha_1^N}, \cdots, v_{\alpha_{N^2}^N} ]^T$,
which lists the $N^2$ entries by the same rule as in \eqref{xwN2},
and
$ {\cal D}=I_N\otimes {\cal R}_N-2 I_{N^2}.
$
 We can show
 \begin{align}\label{sdevn}
 d{\bf v}(t)= {\cal M}^{tr}_N {\bf v}(t) dt +\sigma {\cal D} d{\bf  w}(t).
 \end{align}
Using \eqref{sdevn} we obtain the following lemma.

\begin{lemma}\label{lem:vBound}
 We have
$
\sup_{k,0\le t\le T}E|{ v}_{\alpha_k^N}(t)|^2\le C\ln N/(N^2N_1).
$
\end{lemma}

\proof
The matrix
${\cal D}^2$ has all its eigenvalues (including algebraic multiplicities) given by the $N$ values
  $(2-\lambda_1)^2, \cdots, (2-\lambda_N)^2  $, each repeated $N$ times (see
  \cite[Theorem 4.2.12]{HJ91}), where $\lambda_i$ is defined by \eqref{lamkN}.
  Hence ${\cal D}^2\le 16I_{N^2}$.
 By ${\bf v}(0)=0$, we have
${\bf v}(t)= \sigma \int_0^t { \Phi^{tr} (t,\tau)} {\cal D}
d{\bf  w}(\tau)$, and
\begin{align}
{\mathbb E}|{\bf v}(t)|^2
&\le (16\sigma^2/N_1) \mbox{Tr} \int_0^t e^{2\int_\tau^t {\cal M}^{tr}_N (r)dr}d\tau\le C  \ln N/N_1,
\end{align}
where we have used \eqref{trlog} to get the last bound.
By  symmetry of the components in ${\bf v}$, the lemma follows.  \qed

\begin{remark} \label{re:verticalv}
The  bound in Lemma \ref{lem:vBound} is still valid after redefining $v_{\alpha_k^N}$
along the vertical direction of the array of nodes on ${\mathbb T}^2_{1/N}$.
\end{remark}

\begin{lemma}\label{lem:NNhxbx}
We have
\begin{align}
&\sup_{t,k}\EE\Big|N^2 \sum_{k':\alpha_{k'}^N\sim  \alpha_{k}^N   }(\hat x^{(N_1)}_{\alpha_{k'}^N}
- \hat x^{(N_1)}_{\alpha_k^N} )- N^2 \sum_{k': \alpha_{k'}^N \sim \alpha_{k}^N  }(\bar x_{\alpha_{k'}^N}
-\bar x_{\alpha_k^N} )\Big|^2  %\\&\qquad
 \le C N^2\ln N/N_1. \notag
\end{align}
\end{lemma}

\begin{proof}
The lemma follows from
 Lemma \ref{lem:vBound} and Remark \ref{re:verticalv}.
\end{proof}

 Denote
$
\xi_{\alpha_k^N}(t)= N^2 \sum_{k': \alpha_{k'}^N \sim \alpha_{k}^N  }(\bar x_{\alpha_{k'}^N}
-\bar x_{\alpha_k^N} )(t)
$.
We will further check the asymptotic behavior of $\xi_{\alpha_k^N}(t)$. For convenience of analysis,  we split the summation in  $\xi_{\alpha_k^N}$ into two parts corresponding to neighbors along  two directions of the array of nodes on ${\mathbb T}^2_{1/N}$.  The estimate in Lemma \ref{lem:delN} below will mean that
equation \eqref{checkxd}  may be regarded as  a space discretized version of the  parabolic equation of $m(t,\alpha)$. Along the horizontal direction of the array for ${\mathbb T}_{1/N}^2$,  define
\begin{align}
\theta_{\alpha_k^N}(t):= N^2 (\bar x_{\alpha_{k_+}^N}
+\bar x_{\alpha_{k_-}^N}
-2\bar x_{\alpha_k^N} )(t). \label{theant}
\end{align}
For $\alpha\in [0,1)^2\subset {\mathbb R}^2$,  denote $\alpha=(\alpha^{(1)}, \alpha^{(2)})$.
Write the second order partial derivatives $\partial_{\alpha^{(1)}}^2 m(t,\alpha)$ and $\partial_{\alpha^{(2)}}^2 m(t,\alpha)$. We proceed to  show that $\theta_{\alpha_k^N}(t) $  approaches $\partial_{\alpha^{(1)}}^2 m(t,\alpha_k^N) $ as $N\to \infty$.

\begin{lemma}
Let $\theta_{\alpha_k^N}(t)$, $1\le k\le N^2$, be defined by \eqref{theant}. Then $(\theta_{\alpha_1^N}, \cdots, \theta_{\alpha_{N^2}^N }  ) $ is the unique solution to the ODE system:
\begin{align}
\frac{d}{dt}\theta_{\alpha_k^N}(t) =\ & [A-B^2R^{-1} P(t)+D] \theta_{\alpha_k^N} (t)\notag \\
&
- B^2R^{-1}N^2[ S(t,\alpha_{k_+}^N)+ S(t,\alpha_{k_-}^N) -2S(t,{\alpha_k^N})  ]  \notag    \\
&+ \kappa N^2 \sum_{k': \alpha_{k'}^N\sim \alpha_k^N}( \theta_{\alpha_{k'}^N}  - \theta_{\alpha_k^N}  )(t) ,\quad
1\le k\le N^2,    \label{checkxxd}
  \end{align}
  with initial condition $\theta_{\alpha_k^N}(0)=N^2 [m_0({\alpha_{k_+}^N})
+m_0({\alpha_{k_-}^N})
-2m_0({\alpha_k^N} )]$.
 \end{lemma}

 \begin{proof}
 We take
 $\alpha_k^N, \alpha_{k_+}^N$ and $ \alpha_{k_-}^N  $ in \eqref{checkxd}
 to obtain three ODEs, which are multiplied by $N^2$ and combined together to derive \eqref{checkxxd}. The solution to \eqref{checkxxd} is unique since it is a linear ODE system with bounded coefficients  when $S$ is given.
 \end{proof}

\begin{lemma} \label{lem:delN}
We have
\begin{align*}
\lim_{N\to \infty}\sup_{t\le T,k\le N^2}\big|\theta_{\alpha_k^N}(t) - \partial_{\alpha^{(1)}}^2 m(t,\alpha)
|_{\alpha=\alpha_k^N}\big|=0.
\end{align*}
\end{lemma}

\begin{proof}
To analyze \eqref{checkxxd}, we introduce the parabolic equation
\begin{align}
\partial_t \partial_{\alpha^{(1)}}^2  m(t,\alpha)= & [A-B^2R^{-1} P(t)+D]
\partial_{\alpha^{(1)}}^2  m(t,\alpha )
\notag \\
&
- B^2R^{-1}  \partial_{\alpha^{(1)}}^2  S(t,\alpha)  dt +
\kappa \Delta (\partial_{\alpha^{(1)}}^2  m(t,\alpha) ),
 \label{pdem2}
\end{align}
where the initial condition is $\partial^2_{\alpha^{(1)}}m(0,\alpha)$ (periodic in $\alpha$).  We view $\partial_{\alpha^{(1)}}^2 m(t,\alpha)$ as the solution of a second order parabolic equation obtained by differentiating \eqref{pde4m} twice.
There is a unique solution to the resulting linear parabolic equation subject to the  initial condition $\partial^2_{\alpha^{(1)}}m(0,\alpha)$.   Denoting
$\chi(t,\alpha)=\chi_\alpha(t)=\partial_{\alpha^{(1)}}^2m(t,\alpha)$ and taking $\alpha=\alpha_k^N$ in \eqref{pdem2}, we obtain  an ODE system with perturbation terms:
\begin{align}
\frac{d}{dt} \chi_{\alpha_k^N}(t) =& [A-B^2R^{-1} P(t)+D]\chi_{\alpha_k^N}(t)
\label{dchik} \\
&- B^2 R^{-1} N^2 [ S(t,{\alpha_{k_+}^N})+ S(t,{\alpha_{k_-}^N}) -2S(t,{\alpha_k^N})  ] \notag\\
 &+  \kappa  N^2 \sum_{k': \alpha_{k'}^N\sim \alpha_k^N}( \chi_{\alpha_{k'}^N}   - \chi_{\alpha_k^N}  )(t)  +\epsilon^o_{\alpha_k^N}(t)  ,\qquad 1\le k\le N^2, \notag
\end{align}
where
\begin{align*}
\epsilon^o_{\alpha_k^N}(t)&=-  B^2R^{-1}  \partial_{\alpha^{(1)}}^2  S(t,\alpha)|_{\alpha=\alpha_k^N}\\
&\quad + B^2 R^{-1} N^2 [ S(t,{\alpha_{k_+}^N})+ S(t,{\alpha_{k_-}^N}) -2S(t,{\alpha_k^N})  ]\notag\\
 &\quad+ \kappa \Delta (\partial_{\alpha^{(1)}}^2  m(t,\alpha) )|_{\alpha=\alpha_k^N}- \kappa  N^2 \sum_{k': \alpha_{k'}^N\sim \alpha_k^N}( \chi_{\alpha_{k'}^N}   - \chi_{\alpha_k^N}  )(t).\notag
\end{align*}
We proceed to compare the solutions of \eqref{checkxxd}
and \eqref{dchik}.
It is clear that  $\epsilon^o_{\alpha_k^N}(t)$ is continuously differentiable on $[0,T]$ and  $\sup_{t,k}| \epsilon^o_{\alpha_k^N}(t)|=o(1) $ as  $N\to \infty$, which is due to the approximation of
$\partial_{\alpha^{(1)}}^2  S(t,\alpha )$ and $\Delta( \partial_{\alpha^{(1)}}^2 \chi(t,\alpha)) $
 by the second order difference quotients under (H3).
For \eqref{checkxxd} and \eqref{dchik},  denote
  \begin{align}
\epsilon_{\alpha_k^N}(t)=\theta_{\alpha_k^N}(t)- \chi_{\alpha_k^N} (t),
\qquad  \epsilon_N= \sup_{k,t}| \epsilon_{\alpha_k^N}(t) |.   \notag
 \end{align}
Then   we have
\begin{align}
\frac{d}{dt}\epsilon_{\alpha_k^N} = & [A-B^2R^{-1} P(t)+D] \epsilon_{\alpha_k^N} \notag \\
&+ \kappa  N^2 \sum_{k': \alpha_{k'}^N\sim \alpha_k^N}  ( \epsilon_{\alpha_{k'}^N}   - \epsilon_{\alpha_k^N}  )   -\epsilon^o_{\alpha_k^N},
 \quad 1\le k\le N^2,  \notag
\end{align}
where $\epsilon_{\alpha_k^N}(0)= N^2 [m_0({\alpha_{k_+}^N})
+m_0({\alpha_{k_-}^N})
-2m_0({\alpha_k^N} )] -\partial^2_{\alpha^{(1)}} \widetilde m_0(\alpha_k^N)$.

 Thus following the proof of Lemma \ref{lem:xbBound},  we may construct two auxiliary ODEs by setting $c_1=\pm \sup_{k,t}|\epsilon^o_{\alpha_k^N}(t)|$ and initial condition $\pm \sup_k|\epsilon_{\alpha_k^N}(0)|  $ in
\eqref{yc1} to bound  $\epsilon_{\alpha_k^N}(t)$ both from above and from below on $[0,T]$.
So we have $\epsilon_N=o(1)$ in view of the estimate
$\sup_{k,t}(| \epsilon^o_{\alpha_k^N}(t)| +|\epsilon_{\alpha_k^N}(0)|)  =o(1) $ under (H1) and (H3).
\end{proof}

\begin{remark}\label{re:sitm}
If we define $\theta'_{\alpha_k^N}(t)$    using two neighbors along the vertical direction instead of \eqref{theant}, then we similarly obtain
$$
\lim_{N\to \infty}\sup_{ t\le T,  k\le N^2}|\theta'_{\alpha_k^N}(t)-
\partial_{\alpha^{(2)}}^2m(t,\alpha)|_{\alpha=\alpha_k^N}|
=0.
$$
\end{remark}

\begin{remark}\label{re:Tm}
For the ring network ${\mathbb T}_{1/N}$, Lemma \ref{lem:NNhxbx} holds with an upper bound $CN^3/N_1$, and Lemma
\ref{lem:delN} still holds after replacing
$\partial^2_{\alpha^{(1)}} m(t,\alpha)$ by $\partial_\alpha^2 m(t,\alpha)$, and $k\le N^2$ by $k\le N$.
\end{remark}

\subsection{The proof of Theorem \ref{prop:Jbound}}
\label{sec:sub:Jiproof}

\proof We use
Lemmas \ref{lem:NNhxbx},  \ref{lem:delN}, Remarks \ref{re:sitm}, \ref{re:Tm}, and (H4) to  obtain
\begin{align}\label{ssupxhxh}
\sup_{N,N_1}\sup_{t,k} \EE\Big| N^2 \sum_{\alpha_{k'}^N\sim \alpha_k^N} (\hat x^{(N_1)}_{\alpha_{k'}^N} - \hat x_{\alpha_k^N}^{(N_1)} )(t)\Big|^2\le C
\end{align}
for some constant $C$, which combined with Lemma \ref{prop:exhBound} yields
$\sup_{t,k,i} \EE |\hat x^i_{\alpha_k^N}(t)|^2= O(1)$.
Then we obtain the cost bound $\hat C_0$  in Theorem
\ref{prop:Jbound}.
\qed

\section{The $\epsilon$-Nash equilibrium}
\label{sec:eNE}

This section considers the network cases of  both ${\mathbb T}_{1/N}$ and ${\mathbb T}_{1/N}^2$.
To establish the $\epsilon$-Nash equilibrium property,
consider  alternative strategies for a single agent.
Let the space ${\cal U}^c_{\rm ind }$ consist of individual feedback strategies $\varphi(t,{\bf x})$, $0\le t\le T$, with ${\bf x}$ being the state vector of the overall population (allowing centralized state information), such that $\varphi$ is continuous in $(t, {\bf x})$ and Lipschitz continuous in ${\bf x}$.
In particular,  each decentralized strategy $\hat u^i_{\alpha_k^N}$ belongs to ${\cal U}^c_{\rm ind }$.
If  $ u^l_{\alpha_j^N}\in {\cal U}^c_{\rm ind  } $ for all $j,l$, the set of individual strategies $ ( u^i_{\alpha_k^N}, u^{-i}_{\alpha_k^N},  {\bf u}_{-\alpha_k^N}   )  $ ensures a well defined solution to \eqref{dxat3}, where $ u^{-i}_{\alpha_k^N}$ and $ {\bf u}_{-\alpha_k^N}  $ follow the notation  in Section \ref{sec:lq}.
Define the strategy set  $ {\cal U}_{ {\alpha_k^N}}^{i,\hat C_0} $
consisting of all $u_{\alpha_k^N}^i \in {\cal U}^c_{\rm ind  }  $ such that   $J_{\alpha_k^N}^i(u^i_{\alpha_k^N},\hat u^{-i}_{\alpha_k^N},  \hat{\bf u}_{-\alpha_k^N}  )$ satisfies condition
\eqref{pcbc0}.

For the performance estimate, we need to approximate
the two terms
$$
x_{\alpha_k^N}^{(N_1)}(t), \quad  N^2 \sum_{k':\alpha_{k'}^N\sim  \alpha_{k}^N   }( x^{(N_1)}_{\alpha_{k'}^N}
-  x^{(N_1)}_{\alpha_k^N} )(t) \quad \mbox{in System (B)},
$$
respectively, by
$$
\hat x_{\alpha_k^N}^{(N_1)}(t), \quad  N^2 \sum_{k':\alpha_{k'}^N\sim  \alpha_{k}^N   }(\hat x^{(N_1)}_{\alpha_{k'}^N}
- \hat x^{(N_1)}_{\alpha_k^N} )(t) \quad \mbox{in System (A)}.
$$
Denote
\begin{align} \label{zetkn}
\zeta_{\alpha_k^N}(t) =  N^2 \sum_{k':\alpha_{k'}^N\sim  \alpha_{k}^N   }( x^{(N_1)}_{\alpha_{k'}^N}
-  x^{(N_1)}_{\alpha_k^N} )(t)  - N^2 \sum_{k':\alpha_{k'}^N\sim  \alpha_{k}^N   }(\hat x^{(N_1)}_{\alpha_{k'}^N}
- \hat x^{(N_1)}_{\alpha_k^N} )(t) .
\end{align}
In the following, without loss of generality, we consider an alternative strategy for agent ${\cal A}^1_{\alpha_1^N}$ so that $( u^1_{\alpha_1^N}, \hat u^{-1}_{\alpha_1^N},  \hat {\bf u}_{-\alpha_1^N}   )$
is applied by the population in System (B).

\begin{lemma}\label{lem:xzeta}
For System (B) with $( u^1_{\alpha_1^N}, \hat u^{-1}_{\alpha_1^N},  \hat {\bf u}_{-\alpha_1^N}   )$ and $u^1_{\alpha_1^N}\in
{\cal U}^{1,\hat C_0}_{{\alpha_1^N}} $, we have
\begin{align}
&\sup_{u_{\alpha_1^N}^1 \in {\cal U}_{ {\alpha_1^N}}^{1,\hat C_0}  }
\sup_{t\le T} \EE |x_{\alpha_1^N}^{(N_1)}(t)- \hat x_{\alpha_1^N}^{(N_1)}(t)|^2= O(1/N_1^2),
\label{xxn1}\\
&\sup_{u_{\alpha_1^N}^1\in  {\cal U}_{{\alpha_1^N}}^{1,\hat C_0}  }  \sup_{t\le T} \EE |\zeta_{\alpha_1^N}(t)|^2 = O(N^4/N_1^2).\label{chn1}
\end{align}
\end{lemma}

\begin{proof}
See Appendix A.
\end{proof}

By Lemmas  \ref{lem:xbm5}, \ref{lem:xchBound},   \ref{lem:NNhxbx},  \ref{lem:delN},    and \ref{lem:xzeta}, Remarks \ref{re:sitm} and \ref{re:Tm},
we further obtain the following corollary.

\begin{corollary}\label{cor:xm}
For  the cases of both  ${\mathbb T}_{1/N}$ and ${\mathbb T}_{1/N}^2$,
as $N\to \infty$, we have
\begin{align}
&\sup_{u^1_{\alpha_1^N}\in  {\cal U}^{1,\hat C_0}_{{\alpha_1^N}}}\sup_{t} \EE |x^{(N_1)}_{\alpha_1^N}(t)- m(t,\alpha_1^N)|^2=o(1),\\
&\sup_{u^1_{\alpha_1^N}\in {\cal U}^{1,\hat C_0}_{{\alpha_1^N}}} \sup_{t} \EE \Big|  N^2 \sum_{k':\alpha_{k'}^N\sim  \alpha_{1}^N   }( x^{(N_1)}_{\alpha_{k'}^N}
-  x^{(N_1)}_{\alpha_1^N} )-\Delta m(t,\alpha_1^N)\Big|^2 =o(1).
\end{align}
%as $N\to \infty$.
\end{corollary}

\begin{theorem} Under assumptions (H1), (H2), (H3), and (H4),
the set of decentralized strategies
$\{\hat u^i_{\alpha_k^N}, 1\le i\le N_1, 1\le k\le N^*\}$ in \eqref{uhati}  is an  $\epsilon_N$-Nash equilibrium for the finite population of $N_1N^*$ agents, i.e.,
\begin{align}
 J_{\alpha_k^N}^i(\hat u^i_{\alpha_k^N},\hat u^{-i}_{\alpha_k^N}, \hat {\bf u}_{-\alpha_{k}^N}  )\le  \inf_{u^i_{\alpha_k^N}\in {\cal U}^c_{\rm ind  }}J_{\alpha_k^N}^i(u^i_{\alpha_k^N},\hat u^{-i}_{\alpha_k^N}, \hat {\bf u}_{-\alpha_{k}^N}  ) +\epsilon_N, \quad \forall i,k, \notag
\end{align}
 where
$\epsilon_N =o(1)$ as $N\to \infty$.
\end{theorem}

{\it Proof}.
Without loss of generality, consider agent ${\cal A}_{\alpha_1^N}^1$ taking  an alternative strategy $ u^1_{\alpha_1^N}$. By Theorem \ref{prop:Jbound}, it suffices to consider
$u^1_{\alpha_1^N} \in {\cal U}^{1,\hat C_0}_{{\alpha_1^N} } $.
We rewrite the equation of ${\cal A}_{\alpha_1^N}^1$ in System (B)
with $n=1$:
\begin{align}\label{xsysB}
dx^1_{\alpha_1^N} =\ &[ Ax^1_{\alpha_1^N}
+ Bu^1_{\alpha_1^N}
+ D  x_{\alpha_1^N}^{(N_1)}]dt \\
&+
\kappa N^2 \sum_{k:\alpha_k^N\sim \alpha_1^N } (x_{\alpha_k^N}^{(N_1)}
-x_{\alpha_1^N}^{(N_1)})dt +\sigma dw_{\alpha_1^N}^1.\notag
\end{align}

To estimate $J^1_{\alpha_k^N}$,
we introduce an auxiliary  optimal control problem with dynamics and cost:
\begin{align}\label{ztab}
dz(t)&=[Az(t)+Bu^\dag(t) +Dm (t,\alpha_1^N)]dt + \kappa \Delta m(t,{\alpha_1^N}) dt  +\sigma dw^1_{\alpha_1^N} (t),\\
J^*(u^\dag)&=  \EE \int_0^T\Big[ \|z(t)-\Gamma_1 m(t,\alpha_1^N)-\Gamma_2\eta\|_Q^2   + \|u^\dag(t)\|_R^2 \Big]dt \notag \\
&\qquad\quad + \EE  \|z(T)-\Gamma_{1f} m(T,\alpha_1^N)-\Gamma_{2f}\eta\|_{Q_f}^2,
\notag
\end{align}
where $z(0)= x_{\alpha_1^N}^1(0)$.
 By the solution of the best response control problem in \eqref{lqdyn}--\eqref{lqcos},
 the optimal control law in \eqref{ztab} is given by
\begin{align*}
\hat u^\dag(t)= - B R^{-1} [ P(t)z(t)  + S(t,\alpha_1^N)].
\end{align*}

For the set of decentralized strategies in \eqref{uhati}, we use
Corollary \ref{cor:xm} and elementary $L^2$ error estimates  to show
\begin{align}\label{jajso1}
|J_{\alpha_1^N}^1(\hat u^1_{\alpha_1^N},\hat u^{-1}_{\alpha_1^N}, \hat {\bf u}_{-\alpha_{1}^N}  ) - J^*(\hat u^\dag) |=o(1), \quad \mbox{as}\ N\to \infty.
\end{align}
Throughout the proof, each $o(1)$ term may be estimated without depending on which agent has been selected for an alternative control.
When $\hat u^1_{\alpha_1^N}$ is replaced by a strategy $u^1_{\alpha_1^N}
\in  {\cal U}^{1,\hat C_0}_{ {\alpha_1^N} }   $
to generate  System (B), we need to show that the cost may be reduced by at most $o(1)$.
We write
\begin{align*}
 &J_{\alpha_1^N}^1(u^1_{\alpha_1^N},\hat u^{-1}_{\alpha_1^N}, \hat {\bf u}_{-\alpha_{1}^N}  )\\
 & = \EE \int_0^T\Big( \|x_{\alpha_1^N}^1(t)-\Gamma_1 x^{(N_1)}_{\alpha_1^N}(t)-\Gamma_2\eta\|_Q^2   + \|u_{\alpha_1^N}^1(t)\|_R^2\Big )dt \notag \\
&\qquad\quad + \EE  \|x_{\alpha_1^N}^1(T)-\Gamma_{1f} x^{(N_1)}_{\alpha_1^N}(T)-\Gamma_{2f}\eta\|_{Q_f}^2.
\notag
\end{align*}

We write
\begin{align}
 &\EE \int_0^T \|x_{\alpha_1^N}^1(t)-\Gamma_1 x^{(N_1)}_{\alpha_1^N}
 (t)-\Gamma_2\eta\|_Q^2dt
 \label{emz7}  \\
&= \EE \int_0^T \|z(t)-\Gamma_1 m(t,\alpha_1^N) -\Gamma_2\eta\|_Q^2dt \notag \\
&\quad+\EE \int_0^T \|z(t)- x_{\alpha_1^N}^1(t)+\Gamma_1 [  x^{(N_1)}_{\alpha_1^N}(t)
-m(t,{\alpha_1^N})
 ]\|^2_Q dt \notag \\
 &\quad -2 \EE \int_0^T [z(t)-\Gamma_1 m(t,\alpha_1^N) -\Gamma_2\eta ] Q  \notag \\
&\qquad \qquad  \cdot [z(t)- x_{\alpha_1^N}^1(t)+\Gamma_1 (  x^{(N_1)}_{\alpha_1^N}(t)-m(t,{\alpha_1^N})   ] dt. \notag
\end{align}
When System (B) has generated  $u_{\alpha_1^N}^1(t, \omega)$ as a random process, we set $u^\dag(t,\omega)=u_{\alpha_1^N}^1(t, \omega) $ in \eqref{ztab}.
By comparing  $z(t)$ in \eqref{ztab} and $ x_{\alpha_1^N}^1(t) $ in \eqref{xsysB} under the same control process  $u_{\alpha_1^N}^1(t,\omega) $ on $[0,T]$, we use
Corollary \ref{cor:xm} to obtain  %$XXXX {\cal U}_{\alpha_1^N}^{1,C_0} $
\begin{align}
\sup_{u^1_{\alpha_1^N}\in {\cal U}^{1,\hat C_0}_{ {\alpha_1^N} }  }\sup_{0\le t\le T} \EE
|z(t)- x_{\alpha_1^N}^1(t)|^2=o(1), \label{supzx7}
\end{align}
as $N\to \infty$.  By \eqref{emz7}, \eqref{supzx7} and Corollary \ref{cor:xm},
we have
\begin{align}
&\sup_{u^1_{\alpha_1^N}\in {\cal U}^{1,\hat C_0}_{{\alpha_1^N} }  }\Big|\EE \int_0^T \|x_{\alpha_1^N}^1(t)-\Gamma_1 x^{(N_1)}_{\alpha_1^N}(t)-\Gamma_2\eta\|_Q^2dt
\notag \\
&
\qquad\qquad\qquad\qquad - \EE \int_0^T \|z(t)-\Gamma_1 m(t,\alpha_1^N) -\Gamma_2\eta\|_Q^2dt\Big|=o(1).\notag
\end{align}
We similarly have
\begin{align*}
&\sup_{u^1_{\alpha_1^N}\in {\cal U}^{1,\hat C_0}_{ {\alpha_1^N} }  }\Big|\EE  \|x_{\alpha_1^N}^1(T)-\Gamma_{1f} x^{(N_1)}_{\alpha_1^N}(T)-\Gamma_{2f}\eta\|_{Q_f}^2-\EE \|z(T)-\Gamma_{1f} m(T,\alpha_1^N) -\Gamma_{2f}\eta\|_{Q_f}^2\Big|\\
&\qquad \qquad =o(1).
\end{align*}
Therefore, we have
\begin{align}
J_{\alpha_1^N}^1(u^1_{\alpha_1^N},\hat u^{-1}_{\alpha_1^N}, \hat {\bf u}_{-\alpha_{1}^N}  )&\ge \EE \int_0^T\Big[ \|z(t)-\Gamma_1 m(t,\alpha_1^N) -\Gamma_2\eta\|_Q^2 +\|u_{\alpha_1^N}^1(t)\|_R^2 \Big] dt \notag \\
 &\quad + \EE \|z(T)-\Gamma_{1f} m(T,\alpha_1^N) -\Gamma_{2f}\eta\|_{Q_f}^2-o(1) \notag \\
 &\ge J^*(\hat u^\dag)-  o(1), \label{jjsta}
\end{align}
where $z(t)$ is generated by \eqref{ztab} using the process $u_{\alpha_1^N}^1(t,\omega) $.
By \eqref{jajso1} and \eqref{jjsta}, it follows that
$$
J_{\alpha_1^N}^1(u^1_{\alpha_1^N},\hat u^{-1}_{\alpha_1^N}, \hat {\bf u}_{-\alpha_{1}^N}  )\ge J_{\alpha_1^N}^1(\hat u^1_{\alpha_1^N},\hat u^{-1}_{\alpha_1^N}, \hat {\bf u}_{-\alpha_{1}^N}  ) -o(1).
$$
This completes the proof.
\qed

\section{Conclusion}

Future work will address  general  existence and uniqueness analyses for nonlinear Laplexion mean field game equation  systems.

\section*{Appendix A}
\renewcommand{\theequation}{A.\arabic{equation}}
\setcounter{equation}{0}
\renewcommand{\thetheorem}{A.\arabic{theorem}}
\setcounter{theorem}{0}
\renewcommand{\thesubsection}{A.\arabic{subsection}}
\setcounter{subsection}{0}

\noindent
{\it Proof of Lemma \ref{lem:SK0}}.
  We  first show that $\Lambda_2$ has its image in $  C_{prd}^{\delta/2, \delta}({\mathbf Q})  $.
 Take an arbitrary function $m\in C_{prd}^{1+\delta/2, 2+\delta}({\mathbf Q})
  $ in \eqref{pde4S}, which is solved as an ODE to give $S$ as a function with period 1 in each space variable:
\begin{align}
S(t,\alpha)=& \Lambda_2 (m)(t)\notag \\
=&- \int_t^T\Phi(t,\tau)[ (Q\Gamma_1-P(\tau)D)  m(\tau,\alpha)
- \kappa P(\tau)  \Delta m(\tau,\alpha)  + Q\Gamma_2 \eta  ]d\tau \label{Stal} \\
& - \Phi(t, T)Q_f [\Gamma_{1f} m(T,\alpha )+\Gamma_{2f} \eta ].\notag
\end{align}
To estimate
$|S|_{\delta/2, \delta; {\mathbf Q}}$ using $ [S]_{\delta/2, \delta; {\mathbf Q}^\circ}$ according to Lemma \ref{lem:qholder},
we observe
\begin{align}
|S(t,\alpha)-S(t,{\alpha'})|
&\le
\int_t^T |\Phi(t,\tau)|\cdot | Q \Gamma_1-P(\tau)D| \cdot|m(\tau,\alpha )-m(\tau, {\alpha'})|
d\tau \notag  \\
&\quad +\int_t^T |\Phi(t,\tau)|\cdot\kappa |P(\tau)|\cdot |\Delta m(\tau,\alpha)- \Delta m(\tau, {\alpha'})| d\tau  \notag \\
&\quad + |\Phi(t, T) Q_f \Gamma_{1f}|\cdot |m(T,\alpha) -m(T,{\alpha'})|.\label{ssaa}
\end{align}
By Lemma \ref{lem:qholder} and periodicity of $S$, it suffices to consider $\alpha, \alpha'\in {\mathbb R}^\iota$ with $ d_{max}(\alpha, \alpha'):=\max_{i\le \iota}  |\alpha_i-\alpha_i'|\le  1/2$.
 Hence by \eqref{ssaa}, we have
\begin{align}
& \sup_{\alpha\ne \alpha', t\in [0,T]}
d_{max}(\alpha, \alpha')^{-\delta} {|S(t,\alpha)-S(t,{\alpha'})|} \notag  \\
&\le T\sup_{t,\tau\in [0,T]} |\Phi(t, \tau)|
\cdot( | Q \Gamma_1-PD|_0 +\kappa |P|_0) \notag  \\
&\quad\times \sup_{\tau,\alpha,\alpha'}d_{max}(\alpha, \alpha')^{-\delta}
\Big(|m(\tau,\alpha)-m(\tau, {\alpha'})| + |\Delta m(\tau, \alpha)-\Delta m(\tau, \alpha')|\Big)\notag  \\
&\quad+\sup_t |\Phi (t, T) Q_f \Gamma_{1f} | \cdot
\sup_{\alpha,\alpha'}d_{max}(\alpha, \alpha')^{-\delta} |m(T,\alpha)-m(T,{\alpha'})|\notag    \\
&\le T\Big\{ \sup_{t,\tau} |\Phi(t, \tau)|
\cdot( | Q \Gamma_1-PD|_0 +\kappa |P|_0) \Big\}\notag \\
&\quad \times \Big( \sum_i |\partial_{\alpha_i} m|_{0;{\bf Q}}  +
\sup_{\tau,\alpha,\alpha'}d_{max}(\alpha, \alpha')^{-\delta} |\Delta m(\tau, \alpha) -\Delta m(\tau, \alpha') |\Big)\notag \\
&\quad  + T\sup_t |\Phi (t, T) Q_f \Gamma_{1f} | \sup_{\tau,\alpha,\alpha'}d_{max}(\alpha, \alpha')^{-\delta}  |\partial_\tau m(\tau, \alpha)- \partial_\tau m(\tau, \alpha')  |\notag \\
&\quad + \sup_t |\Phi (t, T) Q_f \Gamma_{1f}| \cdot \sum_i |\partial_{\alpha_i} m(0,\cdot)|_{0;{\mathbb R}^\iota}.
\label{ssaad}
\end{align}
Next for $0\le t<t'\le T$, we have
\begin{align*}
S(t,\alpha)-S(t',\alpha)
=& -\int_t^{t'}\Phi(t,\tau)[ (Q\Gamma_1-P(\tau)D)  m(\tau,\alpha)  -\kappa P(\tau)  \Delta m(\tau,\alpha)]d\tau\\
&+\int_{t'}^T [\Phi(t',\tau)-\Phi(t,\tau) ] \cdot [ (Q\Gamma_1-P(\tau)D)  m(\tau,\alpha)    -  \kappa P(\tau) \Delta m(\tau,\alpha)]d\tau\\
&  + [\Phi(t', T)-\Phi(t,T)] Q_f [\Gamma_{1f}m(T,\alpha )+\Gamma_{2f} \eta ].
\end{align*}
We  use the differentiability  of $\Phi$ to estimate
\begin{align}
& |t-t'|^{-\delta/2}\cdot|S(t,\alpha)-S(t',\alpha)| \notag  \\
&\le
\Big( T^{1-\delta/2}
\sup_{t,\tau }|\Phi(t,\tau)|+ T^{2-\delta/2} \sup_{t,\tau}|\partial_t \Phi(t, \tau)|\Big) \notag  \\
&\qquad \times
(|Q\Gamma_1-PD|_0+\kappa |P|_0) (|m|_{0;{\bf Q}} + |\Delta m|_{0;{\bf Q}})\notag \\
& \quad + T^{1-\delta/2} \sup_t |\partial_t \Phi(t,T)|
\cdot (|Q_f\Gamma_{1f} |\cdot |m|_{0;{\bf Q}} +|Q_f\Gamma_{2f}\eta|). \label{ssttd}
\end{align}
Hence we conclude
$S=\Lambda_2 (m)\in C_{prd}^{\delta/2, \delta}({\mathbf Q})$.

Now it suffices to consider $z=(t, \alpha)$, $z=(t', \alpha')$
such that $d_{max}(\alpha, \alpha')\le 1/2$. We have
\begin{align*}
\frac{|S(z)-S(z')|}{d_p(z , z')^\delta} &\le
\frac{|S(z)- S(t, \alpha') +S(t,\alpha') -S(z')|}
{d_p(z , z')^\delta }\\
& \le d_{max} (\alpha, \alpha')^{-\delta} | S(z)-S(t, \alpha') |
%\\ &\quad
 + |t-t'|^{-\delta/2} |S(t,\alpha')- S(z')|,
\end{align*}
which combined with \eqref{ssaad}, \eqref{ssttd}
and an estimate of $|S|_{0;{\mathbf Q}}$ using  \eqref{Stal}
yields \eqref{sHolder}. \endproof

We give a comparison  lemma
of differential inequalities, which is used in proving Lemma \ref{lem:xbBound}.
For similar comparison results, see \cite{T12}.

\begin{lemma}\label{lem:odecom}
Suppose $c(t)$ is a nonnegative continuous function on $[0,T]$,  $c_0\ge 0$ is a constant,
and the function $z(t)$ is absolutely continuous on
$[0,T]$.

{\rm (i)} If
$\dot{z}(t)\le c(t) z(t) +c_0$
holds almost everywhere on the set $\{t: z(t)> 0\}$, and  $y(t)$ satisfies the ODE $\dot{y}(t) = c(t) y(t)+c_0$ with $y(0)\ge\max\{ z(0),0\}$.
Then $z(t)\le y(t)$ for all $t\in [0,T]$.

{\rm (ii)} If $\dot{z}(t)\ge c(t) z(t) -c_0 $ holds almost everywhere on the set $\{t: z(t)< 0\}$, and  $y(t)$ satisfies the ODE $\dot{y}(t) = c(t) y(t)-c_0$ with $y(0)\le\min\{ z(0),0\}$.
Then $z(t)\ge y(t)$ for all $t\in [0,T]$.
\end{lemma}

\begin{proof}
For case (i), we may assume non-empty $\{t:z(t)>0\}$  since otherwise the conclusion  is obvious. We prove by contradiction and assume there exists $t_1\in (0, T]$ such that
\begin{align}
z(t_1)>y(t_1),\label{zt1yt1}
\end{align}
which implies for some $t_0\in [0,t_1)$, $z(t)\ge y(t)\ge 0$ holds on $[t_0, t_1]$ with $z(t_0)=y(t_0)$. Define $\xi(t)= z(t)-y(t)$. Then
$ \dot\xi(t)\le c(t)\xi(t)  $ holds on $[t_0,t_1]$ so that $\xi(t)\le 0$ on $ [t_0,t_1]$ by Gronwall's lemma. In particular, $\xi(t_1)\le 0$ holds, implying $z(t_1)\le y(t_1)$, which contradicts \eqref{zt1yt1}.
%We have
%$$
%y(t_1)-z(t_1) =\int_{t_0}^{t_1}[ \dot{y}(s) -\dot{z}(s)] ds \ge
%\int_{t_0}^{t_1}c(s)[ y(s)-z(s)] ds \ge 0,
%$$
%leading to  a contradiction.
Therefore we get the desired comparison inequality. The proof of part (ii) is similar and omitted here.
\end{proof}

\begin{lemma}\label{lem:odeF}
Suppose $F(t)$ is a continuous  $\mathbb{R}^{k\times k}$-valued
function defined on $[0,T]$, and $F(t)F(s)=F(s)F(t)$ for all $s,t\in [0,T]$. Then $ \Phi(t,s)= e^{\int_s^t  F(r) dr}$ is the fundamental solution matrix of the linear ODE
$
\dot{x} = F(t)x,
$
i.e., $\frac{\partial}{\partial t}{\Phi(t,s)} = F(t) \Phi(t,s)$ with $\Phi(s,s)=I_k$.
\end{lemma}

\begin{proof} We have $\frac{1}{\epsilon}[\Phi(t+\epsilon,s)-\Phi(t,s)]=
\frac{1}{\epsilon}(e^{\int_t^{t+\epsilon} F(r)dr} -I)  e^{\int_s^t F(r) dr } \to F(t) \Phi(t,s)$ as $\epsilon\to 0$, where the equality holds since $F(t)F(s)=F(s)F(t)$ for all $s,t$. \end{proof}

\noindent
{\it Proof of Lemma \ref{lem:xzeta}}.
For System (B), we write
\begin{align}
d x_{\alpha_1^N}^{(N_1)}   =& [A-B^2R^{-1} P(t)+D]  x^{(N_1)}_{\alpha_1^N} dt
- B^2R^{-1} S(t,{\alpha_1^N}) dt \notag    \\
&+ \kappa N^2 \sum_{k':\alpha_{k'}^N\sim  \alpha_{1}^N   }( x^{(N_1)}_{\alpha_{k'}^N}
-  x^{(N_1)}_{\alpha_1^N} ) dt +(B/N_1)\tilde u^1_{\alpha_1^N}dt +\sigma d w^{(N_1)}_{\alpha_1^N}  , \label{xn11} \\
 d x_{\alpha_k^N}^{(N_1)}   =& [A-B^2R^{-1} P(t)+D]  x^{(N_1)}_{\alpha_k^N} dt
- B^2R^{-1} S(t,{\alpha_k^N}) dt \notag    \\
&+ \kappa N^2 \sum_{k':\alpha_{k'}^N\sim  \alpha_{k}^N   }( x^{(N_1)}_{\alpha_{k'}^N}
-  x^{(N_1)}_{\alpha_k^N} ) dt  +\sigma d w^{(N_1)}_{\alpha_k^N} , \quad 2\le k\le N^*, \label{xn1k}
\end{align}
 where $x_{\alpha_j^N}^{(N_1)}(0) =m_0(\alpha_j^N )  $ for all $1\le j\le N^*$, and
 \begin{align}\label{utu7}
 \tilde u^1_{\alpha_1^N}(t)=  u^1_{\alpha_1^N}(t)
 + BR^{-1}P(t) x^1_{\alpha_1^N}(t) + BR^{-1} S(t, \alpha_1^N),
 \end{align}
 which is interpreted as a random process generated by the closed-loop System
 (B).
 Denote $\tilde z_{\alpha_k^N }(t)=  x_{\alpha_k^N}^{(N_1)}(t)-  \hat x_{\alpha_k^N}^{(N_1)}(t)$,
 $ \tilde {\bf z}= [ \tilde x_{\alpha_1^N } ,\cdots,    \tilde x_{\alpha_{N^*}^N }   ]^T$ and $\tilde {\bf u}= [  \tilde u^1_{\alpha_1^N}, 0,\cdots,   0  ]^T$.

(i) The case of ${\mathbb T}_{1/N}^2$.
For  $u^1_{\alpha_1^N}\in {\cal U}_{{\alpha_1^N}}^{1,\hat C_0} $, we will give some estimates in terms of
$\EE \int_0^T | \tilde u^1_{\alpha_1^N}(t)|^2 dt$, which will  be shown later to be bounded.
By \eqref{xkn4}, \eqref{xn11} and \eqref{xn1k}, and in analogue to \eqref{xcms5},  we have
$$
d\tilde {\bf z}= {\cal M}^{tr}_N \tilde {\bf z}(t) dt + (B/N_1)\tilde {\bf u}(t)dt,\quad \tilde {\bf z}(0)=0,
$$
which gives
\begin{align}\label{ztb7}
\tilde z_{\alpha_k^N }(t)= [ 0, \cdots, 1, \cdots, 0]\int_0^t \Phi^{tr}(t,\tau)(B/N_1) \tilde {\bf u} (\tau) d\tau,
\end{align}where 1 is the $k$-th entry of the row.
Let $(M)_{ij}$ denote the $(i,j)$-th entry of a matrix $M$. By
\eqref{Ph12ts} and an orthogonal transformation we can show
$ \sup_N\sup_{t,s}|(\Phi^{tr}(t,s))_{ij}|\le C$, so that \eqref{ztb7} gives
\begin{align}\label{ztilde7}
\sup_{t,k}\EE |\tilde z_{\alpha_k^N}(t)|^2\le (C /N_1^2)
\EE \int_0^T  |\tilde {u}_{\alpha_1^N}(\tau)|^2 d\tau ,
\end{align}
which, in view of Lemma \ref{prop:exhBound}, implies that
\begin{align}\label{bdxn17}
\sup_{t,k}\EE |x^{(N_1)}_{\alpha_k^N}(t)|^2\le  (C /N_1^2)
\EE \int_0^T  |\tilde {u}_{\alpha_1^N}(\tau)|^2 d\tau + C_1.
 \end{align}

We further estimate  $\zeta_{\alpha_1^N}$ in \eqref{zetkn}.
We have
\begin{align*}
d\tilde z_{\alpha_1^N }&= [A-B^2R^{-1} P(t)+D]  \tilde z_{\alpha_1^N}dt+ \kappa N^2 \sum_{k':\alpha_{k'}^N\sim  \alpha_{1}^N   }(\tilde z_{\alpha_{k'}^N}
- \tilde z_{\alpha_1^N} ) dt\\
&\quad +(B/N_1)\tilde u^1_{\alpha_1^N}dt,\notag \\
d\tilde z_{\alpha_k^N }&= [A-B^2R^{-1} P(t)+D]  \tilde z_{\alpha_k^N}dt+ \kappa N^2 \sum_{k':\alpha_{k'}^N\sim  \alpha_{k}^N   }(\tilde z_{\alpha_{k'}^N}
- \tilde z_{\alpha_k^N} ) dt,\\
& \qquad 2\le k\le N^2.\notag
\end{align*}

Let $\alpha_{i_+}^N$ and $\alpha_{i_-}^N$ be the right and left neighbors of node $\alpha_i^N$, respectively.
Denote
\begin{align*}
\chi_{\alpha_i^N}(t)=  \tilde z_{\alpha_{i_+}^N }(t)   + \tilde z_{\alpha_{i_-}^N }(t)   - 2\tilde z_{\alpha_i^N }(t)      , \quad
\chi=[\chi_{\alpha_1^N}, \cdots, \chi_{\alpha_{N^2}^N}  ]^T.
\end{align*}
Then by the method of getting \eqref{sdevn}, we derive
\begin{align*}
d\chi = {\cal M}_N^{tr} \chi dt + (B/N_1)\tilde {\bf u}^\star dt ,
\quad \chi(0)=0,
\end{align*}
where the entry within $\tilde {\bf u}^\star$ at the row corresponding to  node $\alpha_j^N$ is nonzero if and only if $\alpha_j^N $ is equal to $\alpha_1^N$ or is as its horizontal neighbor. So $\tilde{\bf u}^\star $ contains  a total of  three nonzero entries: $ \tilde u^1_{\alpha_1^N}$,
$ \tilde u^1_{\alpha_1^N}$, and $ -2\tilde u^1_{\alpha_1^N}$.
We have
\begin{align*}
\chi_{\alpha_1^N}(t)=(B/N_1) \int_0^T [1, 0,\cdots, 0]\Phi^{tr}(t,\tau)
\tilde {\bf u}^\star(\tau) d\tau.
\end{align*}
Hence, by the method of establishing \eqref{ztilde7}, we obtain
\begin{align}
\EE|\chi_{\alpha_1^N}(t)|^2
&= (C/N_1^2)\EE\int_0^T  |\tilde {u}_{\alpha_1^N}(\tau)|^2 d\tau .\label{s7nut}
\end{align}
We again get a bound of the form \eqref{s7nut} when two neighbours along the vertical direction are used, so that we get the bound  for $\zeta_{\alpha_1^N}$ in \eqref{zetkn}:
\begin{align}\label{zetatu7}
\sup_t \EE |\zeta_{\alpha_1^N}(t)|^2\le (CN^{4}/N_1^2)\EE \int_0^T  |\tilde {u}_{\alpha_1^N}(\tau)|^2 d\tau,
\end{align}
which combined with \eqref{ssupxhxh} yields
\begin{align}\label{en2t}
\sup_{t,k} \EE \Big|N^2 \sum_{k':\alpha_{k'}^N\sim  \alpha_{k}^N   }( x^{(N_1)}_{\alpha_{k'}^N}
-  x^{(N_1)}_{\alpha_k^N} )(t)\Big|^2\le C+ (CN^{4}/N_1^2)\EE \int_0^T  |\tilde {u}_{\alpha_1^N}(\tau)|^2 d\tau.
\end{align}
Now by \eqref{bdxn17},  \eqref{en2t}, and System (B), we further obtain
\begin{align}\label{supx7}
\sup_t\EE |x_{\alpha_1^N}(t)|^2 \le C +(CN^4/N_1^2)\EE \int_0^T  |\tilde {u}_{\alpha_1^N}(\tau)|^2 d\tau.
\end{align}
 By  \eqref{utu7}, \eqref{supx7} and  $u^1_{\alpha_1^N}\in
{\cal U}^{1,\hat C_0}_{{\alpha_1^N}} $, we derive
\begin{align}
\EE \int_0^T  |\tilde {u}_{\alpha_1^N}(\tau)|^2 d\tau \le C+ (CN^{4}/N_1^2) \EE \int_0^T  |\tilde {u}_{\alpha_1^N}(\tau)|^2 d\tau,\notag
\end{align}
so that under (H4),
$$
\sup_N \sup_{u^1_{\alpha_1^N}\in
{\cal U}^{1,\hat C_0}_{{\alpha_1^N}}}  \EE \int_0^T  |\tilde {u}_{\alpha_1^N}(\tau)|^2 d\tau \le C,
$$
which combined with \eqref{ztilde7} and \eqref{zetatu7} implies \eqref{xxn1}
and \eqref{chn1}.

(ii) The case of ${\mathbb T}_{1/N}$.
The estimate is similar and we omit the detail.
\endproof

\section*{Appendix B}
\renewcommand{\theequation}{B.\arabic{equation}}
\setcounter{equation}{0}
\renewcommand{\thetheorem}{B.\arabic{theorem}}
\setcounter{theorem}{0}
\renewcommand{\thesubsection}{B.\arabic{subsection}}
\setcounter{subsection}{0}

The following matrix illustrates ${\mathcal T}_N$ in \eqref{rksum} with $N=4$:
\begin{align}
{\cal T}_4= \left[
\begin{array}{cccc|cccc|cccc|cccc}
0 & 1 & 0 & 1  & 1 & 0& 0 & 0 &0 &0 & 0 & 0 & 1 & 0 & 0 &0  \\
1 & 0 & 1& 0 & 0 & { 1} & 0 & 0 &0 &0 & 0 & 0 & 0 & 1 & 0 &0  \\
0  &  1 & 0& 1 & 0  &0 & 1 & 0& 0&0 &0 &0 & 0 & 0 &1 &0\\
 1 & 0 & 1 & 0  & 0 & 0& 0 & 1 &0 &0 & 0 & 0 & 0 & 0 & 0 &1  \\
 \hline
1 & 0 & 0& 0 & 0 & 1 & 0 & 1 &1 &0 & 0 & 0 & 0 & 0 & 0 &0  \\
0  &  1 & 0& 0 & 1  &0 & 1 & 0& 0&1 &0 &0 & 0 & 0 &0 &0 \\
0 & 0 & 1 & 0  & 0 & 1& 0 & 1 &0 &0 & 1 & 0 & 0 & 0 & 0 &0  \\
0 & 0 & 0& 1 & 1 & 0 & 1 & 0 &0 &0 & 0 & 1 & 0 & 0 & 0 &0  \\
\hline
0  &  0 & 0& 0 & 1  &0 & 0 & 0& 0&1 &0 &1 & 1 & 0 &0 &0 \\
0 & 0 & 0 & 0  & 0 & 1& 0 & 0 &1 &0 & 1 & 0 & 0 & 1 & 0 &0  \\
0 & 0 & 0& 0 & 0 & 0 & 1 & 0 &0 &1 & 0 & 1 & 0 & 0 & 1 &0  \\
0  &  0 & 0& 0 & 0  &0 & 0 & 1& 1&0 &1 &0 & 0 & 0 &0 &1 \\
\hline
1 & 0 & 0 & 0  & 0 & 0& 0 & 0 &1 &0 & 0 & 0 & 0 & 1 & 0 &1  \\
0 & 1 & 0& 0 & 0 & 0 & 0 & 0 &0 &1 & 0 & 0 & 1 & 0 & 1 &0  \\
0  &  0 & 1& 0 & 0  &0 & 0 & 0& 0&0 &1 &0 & 0 & 1 &0 &1 \\
0  &  0 & 0& 1 & 0  &0 & 0 & 0& 0&0 &0 &1 & 1 & 0 &1 &0
\end{array}\right].\notag
\end{align}

 % \bibliographystyle{plain}
%  \bibliography{peterbib,mybib}
% \bibliographystyle{IEEEtran}

%\bibliography{ifacconf}             % bib file to produce the bibliography
                                                     % with bibtex (preferred)

\end{document}